# ACKNOWLEDGEMENTS

I would like to thank the following sets of people, who are to blame for the existence of this thesis (I bear minimal responsibility):

1. My advisor Rebecca McGraw, as well the other two members of my thesis committee, Jennifer Eli and Aditya Adireja, for pushing me outside of my comfort zone and making my thesis far better.

2. K.P. Mohanan, Tara Mohanan, and Jayasree Subramanian, and other members of ThinQ, for working with me to come up with the ideas in this thesis. Also, for completely destroying my bad ideas and making my other ideas far better.

3. David Taylor and Eric Elert, my fellow graduate students, for being my first sounding boards for ideas, whether they were interested in hearing them or not! Also, for making it explicit when my ideas were terrible.

4. My family and friends back home for still liking me (hopefully?) despite me not picking up any of their phone calls or responding to their texts!

5. Arnav and Dhody for willing to argue/fight with me ever since I remember, and being kind enough to lose each time.

6. My parents, Gauri and Uday, and my grandparents, for the extreme, and often un-called for, levels of concern about my life.

7. Oliver Schirokauer, Jim Walsh, Bob Bosch, Jack Calcut, Sandhya Saini, and Rajyashree Sood, my math teachers at Oberlin and Vasant Valley School, for my interest in mathematics.

8. David Glickenstein, Klaus Lux, and my other math professors at The University of Arizona, for giving me insight into the nature of mathematics.

9. Bertrand Russell, Noam Chomsky, David Hilbert, Albert Einstein, John Dewey, Charles

   Darwin, Richard Feynman, Daniel Kahneman, Ian Stewart, David Tall, and George Polya,

   amongst others, whose ideas have been largely unacknowledged in the body of this thesis, but

   form the basis of the way I think about education.

# TABLE OF CONTENTS

## ABSTRACT

Theory Building has been largely ignored in Mathematics Education, especially at the Middle and High School Levels. This thesis focuses on Assumption Digging, a type of Theory Building similar to what Hilbert undertook in his Grundlagen with Euclidean Geometry. The Euclidean Geometry course in schools is a potentially productive for Assumption Digging since, as has been noted, many courses in Euclidean Geometry contain extraneous postulates and unclear definitions. This focus of this thesis is a teaching experiment with 7th-10th grade students in two schools in Western India, which is an Assumption Digging exercise in Euclidean Geometry. The goal of this thesis is to ascertain the various ways students engage in Assumption Digging, and the role played by the instructor in the Assumption Digging sessions. While most of the student engagement with Assumption Digging was at the level of justifying claims which the instructors questioned, there were a few instances of students questioning arguments and claims made by the instructors. The role the instructors played was crucial in two ways. Firstly, deciding what path the session took and, secondly, engaging with student arguments and skepticism.

# INTRODUCTION

There are various possible reasons for teaching mathematics. One is the practical goal of being able to get jobs in an increasingly competitive economy. Another is the goal of developing citizens who have general awareness, including of mathematical ideas, so that they are less susceptible to fake and misleading propaganda.

Alongside these, there is the goal of helping students gain what Bass (2005) calls 'mathematical sensibility'. Broadly speaking, mathematical sensibility refers to the structures, practices and values that characterize mathematics as a discipline (Weiss, 2009). While inquiring into and clarifying this is an important project, it is outside the scope of this thesis. Moving forward, I will be assuming the value of this goal.

One type of knowledge creation practice in Mathematics is Theory Building. Borrowing from Gowers' (2000) distinction between Theory Building and Problem Solving, I use the notion of Theory Building in mathematics as any activity that aims to better understand, clarify, or create a mathematical theory. Many activities of mathematicians have aspects of both Theory Building and Problem Solving. This will be unpacked further in the Literature Review. While Problem Solving has been an important part of Math Education at least since Polya (Polya, 1963, 1978, 1979), Theory Building is less explored.

This thesis is about a particular form of Theory Building, which I'm calling Assumption Digging. I will be exploring Assumption Digging in the context of Euclidean Geometry at a Middle and High School level. The thesis explores Assumption Digging and the teaching of it, on the basis of a teaching experiment carried out in two schools in Western India, with students, almost all of whom had previously studied Euclidean Geometry.

There are various ways to get at the different aspects of theory building. One way unpack it in terms of the abilities, mindsets, values, and understanding developed through theory building exercises, and which help to become better theory builders. For instance, the ability to create and clarify definitions is required to engage in theory building, along with an understanding of what a definition in mathematics is. The mindsets of precision and rigor, along with a high value placed on precision and rigor are equally necessary.

Another way to classify aspects of theory building in mathematics is in terms of the different motivations that prompt the theory building. To illustrate, the motivation to develop Graph Theory came from Euler attempting to solve a real world problem, that of the Bridges of Konigsberg. The motivation to write the Grundlagen came from Hilbert finding flaws in an existing theory, that of Euclidean Geometry.

While slicing Theory Building into the different abilities, understanding, mindsets, and values is an important task, I will only touch on it in this thesis, and only within the context of Assumption Digging. Slicing it into mathematical motivations will give us a better sense of where Assumption Digging fits in to Theory Building.

**Types of Theory Building**

In order to unpack the types of Theory Building, it would be useful to look back to how mathematicians have historically gone about this process. In the following pages, I will attempt to construct a typology of the mathematical motivations which have gone into building various theories in history. This is not to claim that what follows is a complete typology, but it is meant to be illustrative of the mathematical motivations that go into Theory Building. While this thesis is only about a particular type of activity whose mathematical motivations fall into one of these categories, it is useful to see Assumption Digging in the larger context of Theory Building.

**Theory Building Inspired by the Real World.**

The story of Euler and the Bridges of Konigsberg is well known ("Euler", n.d.). To re-state the problem, Konigsberg, a city in Prussia had 7 bridges connecting various land masses, as shown in Figure 1.

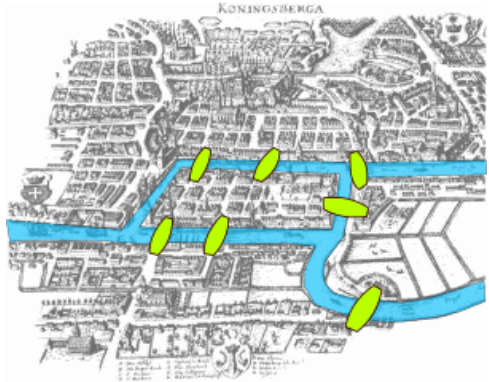

*Figure 1.* Bridges of Konigsberg. From "Wikimedia Commons,"
https://commons.wikimedia.org/wiki/File:Konigsberg_bridges.png. Reprinted with permission under GNU Free Documentation
License.

The task was to find a route through the city which would cross each of the bridges once and only once. Euler showed that such a route is impossible. In order to do that, he represented the city as what we today refer to as a graph, as shown in Figure 2.

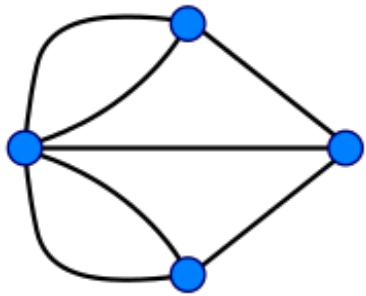

*Figure 2.* Bridges of Konigsberg Graph. From "Wikimedia Commons,"
https://commons.wikimedia.org/wiki/File:K%C3%B6nigsberg_graph.svg. Reprinted with permission under GNU Free
Documentation License.

The land masses were represented by points (vertices/nodes), and the bridges were represented as connections between them (edges/arcs). Once you do that, the answer miraculously appears – if a vertex is neither the start nor end point of the route, then you will have to pass through an even number of edges which are connected to that vertex. Because you have to pass through each edge, all the vertices have an odd number of edges, and there are more than two vertices, no such route exists.

In solving that problem, Euler created Graph Theory, which is even today an important field of study in mathematics, and has application across fields, especially in computation. Even though the theory doesn't refer to bridges and land masses, its inspiration was that problem.

There are many other theories which have directly had inspiration from the real world. For instance, Probability from gambling. An example from the education literature is Paper Folding Geometry (Rao, 1901).

**Generalizing Existing Theories.**

The study of rings is a generalization of the integers. While creating this generalization, we abstract certain properties of the integers such as their closure under addition and multiplication, the existence of the additive identity and additive inverses, the commutativity of addition, and so on. However, we do not require the existence of any sort of ordering relationship between the elements, and do not even require the existence of a multiplicative identity.

Groups generalize various different objects such as the natural numbers, the symmetries of shapes, the symmetries of solutions of polynomials, and so on. Mathematics is littered with examples of theories which emerge from abstracting away aspects of existing theories.

**Changing an Existing Theory.**

Euclid's attempt to prove the parallel postulate was imitated by mathematicians for centuries. Some mathematicians approached the problem in the following manner: if we ignore the parallel postulate, what sorts of consequence do we get? While initially, the goal was to find a contradiction and hence prove the Parallel Postulate, it quickly becomes clear that what results when ignoring the postulate is not necessarily Euclidean Geometry. This resulted in the development of spherical and hyperbolic geometry (Greenberg, 1993).

In spherical geometry, by certain definitions of parallel lines, they do not exist (for instance, the definition as lines which do not meet). By other definitions, every pair of great circles is parallel (for instance, the requirement that there exists a line perpendicular to both). Tweaking the axioms and definitions of existing theories can result in interesting consequences.

**Putting Theories on Firmer Footing.**

Hilbert (1902), and Russel & Whitehead (1912) may have failed in their attempt to create a perfectly sturdy foundation for mathematics. However, the program of delving into foundations and attempting to fix them was a valuable one. A great example is Hilbert's axiomatization of Euclidean Geometry in his Grundlagen. For every claim in each of Euclid's proofs, Hilbert asked the question: why should I accept that? In doing so, Hilbert had to add in a large number of axioms to the ones Euclid had set up (Greenberg, 1993, pp 70-114). Russel did similar things for areas such as Number Theory and Set Theory. In this process of putting theories on firmer footing, we may also realize that some of the existing axioms of the theory are not required since they can be proved from other axioms.

**Assumption Digging**

      I'm using the term Assumption Digging to refer to the sort of activity which Hilbert engaged in with his Grundlagen – finding flaws or gaps in mathematical theories, and fixing them by clarifying and creating definitions and axioms. The motivation here is to place a theory on firmer footing. So, given a conclusion and an argument purporting to prove that conclusion, Assumption Digging is the process of clarifying terms and statements, extracting hidden premises, and clearly spelling out the steps of reasoning to a certain degree of rigor. This degree of rigor could be at the very high level of Russel's Principia (Whitehead & Russell, 1912) or Hilbert's Grundlagen (Hilbert, 2014), or could be at a level which satisfies a skeptical peer. In a classroom, the decision on what constitutes a rigorous enough justification will have to be decided upon. Initially, while students are figuring out what is involved in Assumption Digging, the teacher and the educational materials will probably have to take the lead in this decision. However, as students gain experience, they will hopefully be able to play a larger role.

      **An Analogy.**

      The metaphor of Assumption Digging is related to archaeology. Suppose you are an archaeologist who has found a building buried underground. All you see is the top of the building (the conclusions of the theory). Your goal is to uncover the entire building. As you start digging, you realize that parts of the building are not in great shape. If you wish to dig further, you have to fix those parts and then continue on. Finally, you get to the foundations of the building, which are the axioms, and undefined entities of the theory. Any flaw in the reasoning is the equivalent of a crack in the building which could potentially bring the whole building crashing to the ground.

**Unpacking Assumption Digging.**

To unpack this further, it would be useful to use a diagrammatic representation of Assumption Digging. To make things concrete, I will be using the example which is used in the workshop sessions for which data was collected.

The conclusion (referred to as C) we began with was: Triangles with the same base, and between the same pair of parallel lines, have the same area. The proof for this relies on the following three premises:

Premise 1 (P1): The height of two triangles between the same pair of parallel lines is the same

Premise 2 (P2): The area of a triangle is ½ base x height

Premise 3 (P3): The bases of the two triangles are of the same length

Figure 3 is a diagrammatic representation of this argument.

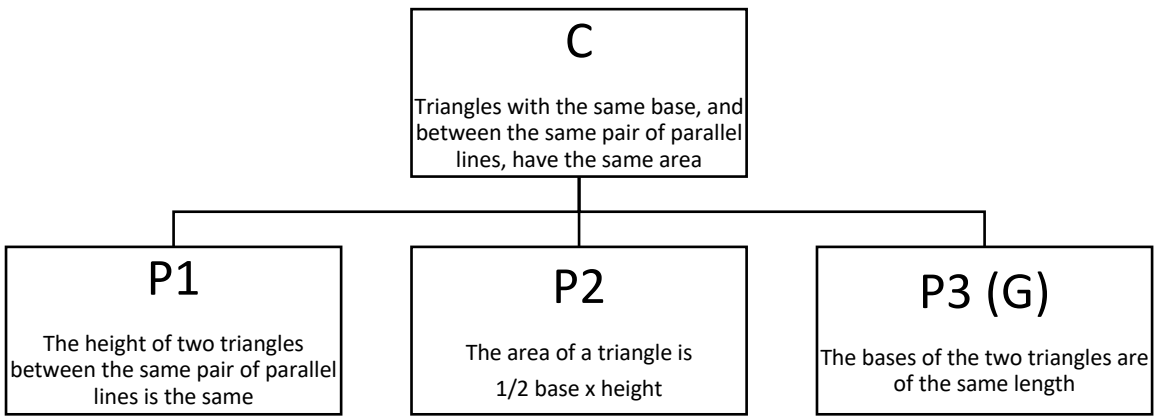

*Figure 3.* Tree representation of the argument above.

At the root of the tree is a conclusion we have begun with. We could begin with more than one conclusion, but the representation would be much messier. Since the particular example used in this thesis starts with just one conclusion, the above representation is sufficient. The proof given for the conclusion (C) is based on the three premises at the second level of the tree.

There are three possible types of premises: Definitions, Axioms, Given, and Claims (which may be theorems – those require a proof).

P3 is given in the statement of the claim we are interested in, and that is indicated by a G.

Taking Premise 2, we could start with the area of a triangle being ½ base x height as the definition of area of a triangle. However, we could also define the area of a triangle in terms of the area of a parallelogram. Here are the premises involved in that argument.

Premise 2-1 (P2-1): Area of a parallelogram is base x height

Premise 2-2 (P2-2): Given a triangle, there exists a parallelogram with the same base and height, and double the area.

Figure 4 is a tree representation of the entire argument so far.

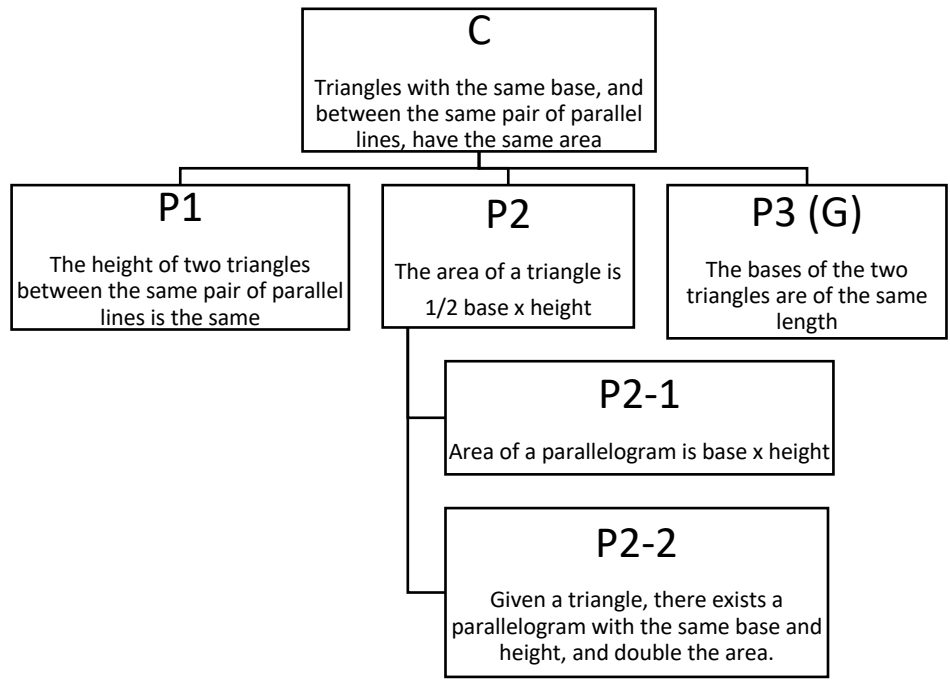

*Figure 4.* Tree representation of the argument above

P2-2 will take some work to prove rigorously – it requires the notion of congruence. With P2-1, we once again have the choice of using that as a definition or we could define the area of a parallelogram in terms of the area of a rectangle.

Taking P1, it will eventually require some form of the parallel postulate, even though we do not get to that in the actual sessions. Figure 5 represents this.

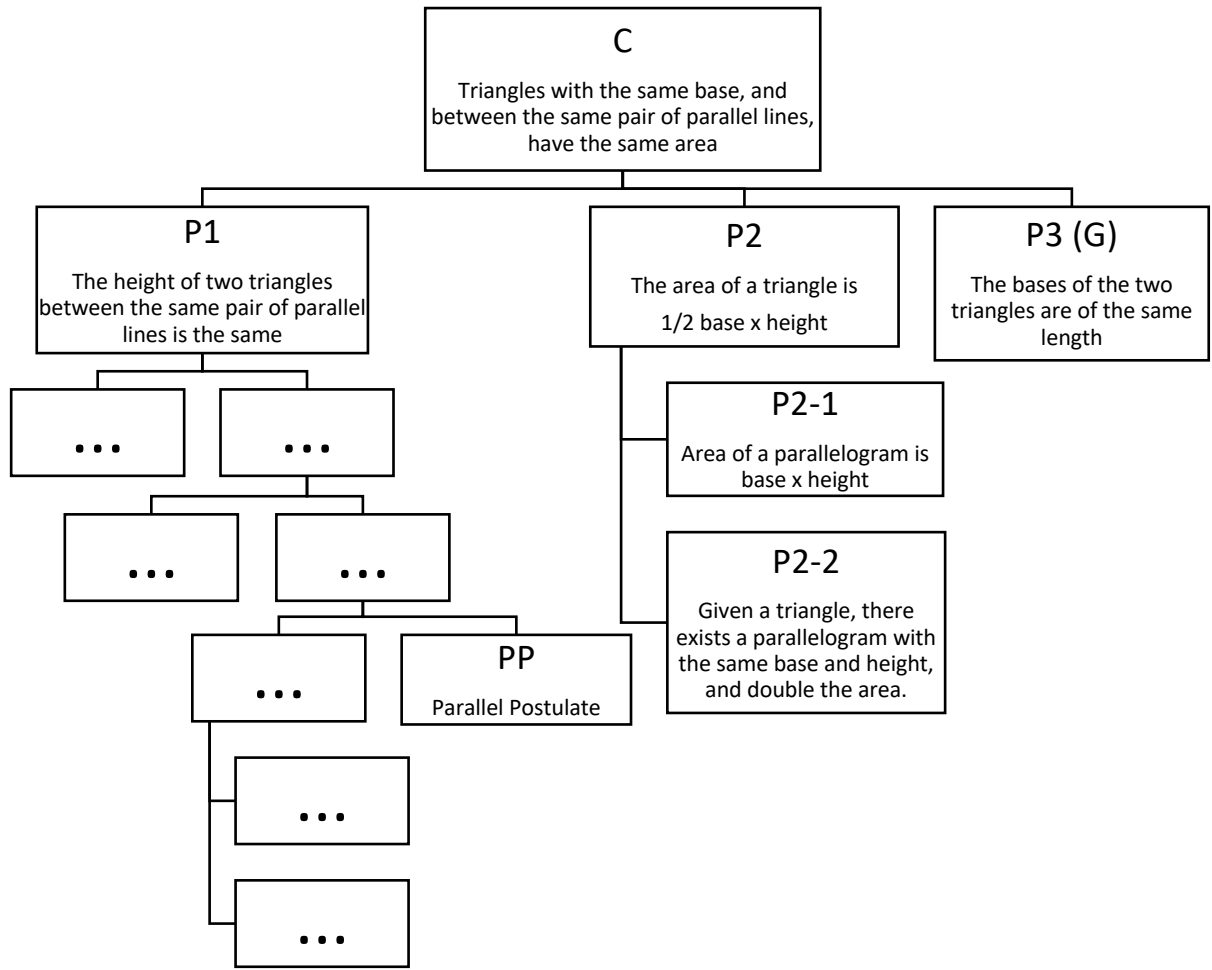

*Figure 5.* Tree representation of an argument going all the way down to the Parallel Postulate

PP represents the parallel postulate in some form. It may well be the case that the same axiom/definition is used at various different points in the diagram – to identify these pictorially would make it extremely messy, and the tree structure would be lost. What is also missing in the diagram is an indication as to what the undefined entities are. While they will be mentioned in the axioms, they also need to be explicitly mentioned. An example of an undefined entity in Euclid's version of his geometry is straight line. It is undefined but is constrained by certain axioms including that given two points (also undefined), there exists exactly one straight line passing through them.

To summarize, the vertices of the tree represent propositions. The root of the tree is a conclusion, while the other vertices are premises on which the conclusion is based. Once the digging is complete, the leaves represent axioms and definitions. The other premises are intermediate conclusions, and could be considered theorems in their own right. The edges of the tree represent the steps of reasoning from each layer of premises to the conclusion one step above.

The creation of this tree, starting with a conclusion and moving downwards towards definitions and axioms, is what is referred to as assumption digging. It may well be the case that your understanding of the claim, and the argument for it from axioms, definitions, and undefined entities, are rigorous. In that case, assumption digging will just be an exercise in explicitly stating that. However, if that is not the case, an important part of assumption digging becomes the search for flaws. These flaws usually come in two flavors, one related to the edges and the other to the vertices:

a.   leaves without justification which have not been stated as axioms or definitions

b.   flaws in the edges – in the steps of reasoning

The assumption digging task, then, becomes an exercise in fixing these flaws. When we see a leaf without justification from a lower layer, there are three possibilities – the leaf contradicts other statements in the theory and should be dropped, the leaf is justified by a lower layer of premises, or the leaf is stated as an axiom or a definition. In order to show the first case, we look for counter-examples which satisfy other claims we wish to accept. It may well be the case that we find a counter-example while attempting a justification.

There is no algorithm to the choice between the other two possibilities. Rather, it would be based on principles such as parsimony, ease of understanding, and beauty. Discussions on such issues could be a useful part of an assumption digging session. However, I will not spend time on them in this thesis since the classroom sessions in this thesis did not get far enough such these considerations would be an issue.

The flaws in the steps of reasoning are easier to get a grasp on. A flaw implies that the conclusion (the higher vertex) doesn't follow from the premises (the lower set of vertices). There are various possibilities this opens up. One possibility is that the particular premises need to be replaced by new ones, or both some edges and their corresponding vertex need to be dropped. Another possibility is that the argument requires extra premises alongside what exists. The other possibility I'm considering is that while the conclusion might follow from the premises, the steps of reasoning are either not clearly articulated or have missing steps. Of course, all of these could be operating together to create further complexity.

### Operationalizing Assumption Digging

For the purposes of this thesis, I will be using the following operational definition of Assumption Digging:

AD1:   Given a claim, extract possible premises on which the claim relies.

AD2:   Attempt to clarify the meanings of the words used in the claim.

AD3:   Continue 1 and 2 recursively on the extracted premises (i.e. for each premise, asking: why should I believe that?)

AD4:   At some point in this recursive exercise, state the premises as axioms and decide that a particular word is an undefined entity with certain properties (which also form some of the axioms of the theory).

In this definition, AD1 and AD2 are unordered. Given this definition, we can now ask

what understanding is required for somebody to be a capable assumption digger. The following

are some of the things a capable assumption digger needs to understand:

a.   Definition (of the iff variety)

b.   Undefined Entity

c.   Axiom

d.   Proof (Mathematical)

In the literature review, I will dig a little deeper into these.

**Classroom Norms for Assumption Digging.**

Without a classroom which encourages productive disciplinary engagement (Engle &

Conant, 2002) amongst students, an assumption digging session would likely not achieve much

of value. Engle & Conant (2002) lay out four principles for creating such a learning

environment:

a.   Problematizing: Encouraging students to take on intellectual problems

b.   Authority: Students are given authority in addressing such problems

c.   Accountability: Students' intellectual work is made accountable to others and to

disciplinary norms

d.   Resources: Students are provided sufficient resources to do all of the above

By its very nature, assumption digging is an exercise in problematizing. The other three norms are important classroom practices which are found in the lessons in this thesis.

**Value of Assumption Digging**

As mentioned earlier in this section, one of the goals often mentioned by mathematics educators is that of getting students to engage in the disciplinary practices of mathematics. This is a part of various curricula. For instance, the Common Core contains 8 standards (Standards for Mathematical Practice, n.d.) for mathematical practice (referred to as MPs). All of them, apart from Mathematical Modeling, have to play some role in an assumption digging session, and three of them have an especially close connection. By its very nature, Assumption Digging involves the construction of viable arguments and critiquing the reasoning of others (MP3). Since we are concerned with making understanding of the justification for a claim more rigorous, there is clear attention to precision (MP6). We are also concerned with structural relationships between statements (MP7). Hence, an Assumption Digging exercise, if successful in what it is trying to achieve, goes a long way in achieving some of the goals of the Common Core (Standards for Mathematical Practice, n.d.).

Similarly, the Indian National Curriculum Framework 2005 position paper on Mathematics (National Council of Educational Research and Training (India), 2006) states that the goals of mathematics education should include students understanding the structure of mathematics, and engaging in meaningful discussions on mathematics. Both of these goals would be partly addressed by a good Assumption Digging session.

Apart from what various educational bodies recommend, one other reason for engaging in assumption digging is that it closely relates to general critical thinking outside of mathematics. When I say critical thinking, I mean the evaluation of claims. While claims outside of mathematics involve forms of reasoning not allowed in mathematics, the strategies involved at a more abstract level are similar.

**Euclidean Geometry**

Although Assumption Digging in Mathematics is something which can happen in any area, Euclidean Geometry is especially fertile for it. The main reason for this is that Euclidean Geometry is often students' first exposure to proving, and often their only experience of proof (Weiss & Herbst, 2015). However, as noted by Weiss and Herbst (2015) regarding Euclidean Geometry in the US, the "geometry course has been criticized for years as a caricature of mathematics, characterized by a proliferation of unnecessary postulates and imprecise definitions, claims that are accepted without proof, and proofs that valorize form over substance" (pp. 225). Whether or not the same is true for students in India would require an analysis of textbooks and more importantly an evaluation of student understanding once they have gone through the course.

# LITERATURE REVIEW

Given that this Thesis is about Assumption Digging in Euclidean Geometry, this Literature Review is in two parts. The first part is concerned with the Literature about Theory Building, Assumption Digging, and related areas. The second is concerned with the teaching of Euclidean Geometry, especially those aspects which lend themselves to Assumption Digging.

## Assumption Digging

### Assumption Digging and RME.

Realistic Mathematics Education (RME) emphasizes giving students 'realistic' situations and problems, which serve as a source for developing concepts, tools and procedures. The word 'realistic' doesn't just refer to the real world. Instead, it refers to problems which are experientially real in students' minds (Van den Heuvel-Panhuizen & Drijvers, 2014, pp. 521-525). This could include problems within formal Mathematics.

In this sense, Assumption Digging, and Theory Building more generally, could be seen within the RME tradition since students are actually engaging with ideas they start with some understanding of rather than being given knowledge by the instructor. Assumption digging could also be seen as a Guided Reinvention activity, where students are re-creating a body of knowledge in an environment guided by the instructor. The role of the instructor in the RME tradition is to be a facilitator, who should use pupils' constructions to help them build new constructions, solve problems, etc. (Wubbels, Korthagen, & Broekman, 1997).

Within this tradition, there has been some work on students building theories such as Dawkin's work on Neutral Axiomatic Geometry (Dawkins, 2015), and Larsen's work on Group Theory (Larsen, 2013). However, most of this work appears to be at an undergraduate level, working with relatively complicated parts of mathematics. However, Zandieh & Rasmussen (2010), which will be described later in this review, could also be considered a Theory Building exercise at the High School level.

**Theory Building.**

Gowers (2000) introduced the distinction between Theory Builders and Problem Solvers as those who see the point of solving problems being to understand the mathematics and its structure better, and those who see the point of understanding mathematics being to solve problems. Gowers sees this distinction as a matter of priorities, and not as an absolute. An example of something clearly involving both problem solving and theory building in about equal measure is the creation of Graph Theory by Euler. It was a problem solving exercise, but required the construction of a theory. He points to areas such as algebraic number theory and geometry as being the abodes of theory builders, where progress "is often a result of clever combinations of a wide range of existing results" (pp. 3). On the other hand, the study of Graph Theory, where the objects can be immediately comprehended, has interesting problems, the solving of which is not that reliant previous work in the field.

While Gowers was talking about the mindset of Mathematicians, Bass (2017), borrowing from Gowers, defines Theory Building as "creative acts of recognizing, articulating, and naming a mathematical concept or construct that is demonstrably common to a variety of apparently different mathematical situations, a concept or construct that, at least for those engaged in the work, might have had no prior conceptual existence." An example he gives of Theory Building is a preschooler being able to abstract the similarity in patterns like the following: ‖**‖**… and <<//<<//<< …, and stating the relationship as something like 'Same, Same, Different, Different'.

Bass' notion of Theory Building captures some aspects of what Gowers says. However, if we follow the consequences of Gowers' statement that Theory Builders are interested in understanding mathematics better, then an example of that would be Hilbert's Grundlagen (Hilbert, 2014). This was Hilbert's attempt at making Euclidean Geometry more rigorous. By finding holes in Euclid's reasoning, he clarified existing definitions and axioms, and added in more axioms where the need for them arose. It isn't clear whether this particular example would fit into Bass' idea of Theory Building, and it definitely doesn't fit into Gower's notion of Problem Solving since there are no new results Hilbert found about Euclidean Geometry. To fix this would either require a clarification of Bass' definition or a new, broader definition of Theory Building. This thesis will not delve into that. Rather, I will be using something similar to Gowers' notion, and say that Theory Building is any activity which aims at better understanding, clarifying, or creating the structure of a mathematical theory. Solving problems in order to better understand the theory is an important part of Theory Building along with things like clarifying definitions, and laying out axioms.

Fawcett (1938) documents possibly the first attempt at the creation of a systematic course aimed at theory building. He starts with high school students listing geometric objects. Over the length of the course, this gets translated to a theory of space with definitions, undefined entities, axioms, and theorems. Fawcett's course is far broader in scope and more ambitious than Assumption Digging. In fact, something like Assumption Digging is happening as a part of his course when the participants extract the structure of the theory from claims they make about the objects.

**Understanding Related to Assumption Digging.**

Looking at the definition of Assumption Digging from the Introduction, we see that it requires an understanding of definitions, axioms, undefined entities, and proofs. I will focus on the literature on Definitions and Proofs. The other aspects appear in those sets of literature.

### *Definitions.*

The distinction between Concept Image and Concept Definition introduced by Tall and Vinner (Tall & Vinner, 1981) is a useful notion in that it points to the need to deal with concepts both as formal entities as well as treating them as the intersection of various properties. The Concept Image is the total cognitive structure associated with a concept. This includes mental pictures, and associated properties and processes. The Concept Definition, on the other hand, is the words used to specify the concept.

Mariotti & Fischbein (1997) discuss the distinction between two types of definitions – those which are the 'basic objects' of the theory and those which are new elements within a theory. The basic objects, which are the same as what I have called undefined entities, have a close relationship with the axioms of the theory. For instance, points and lines in Euclidean Geometry do not have a definition – rather, they are concepts constrained by certain axioms. In the context of Geometry, they are also concerned with the 'figural' and the 'conceptual', where the figural is to do with the relationship between geometrical objects and their counterparts in space, while the conceptual is to do with the abstract nature of geometrical reasoning.

Putting these together, they propose a pedagogy for coming up with definitions via a classroom discussion. This consists of:

- Observing

- Identifying the main characteristics

- Stating properties based on them

- Returning to observations to check

The role of the teacher here is very important, especially in getting students to harmonize the conceptual and figural aspects of the object they are defining. There is a close relationship between the activities in this paper and Assumption Digging since students are coming up with definitions for which they already have a concept.

Another closely related framework is put forward by Zandieh & Rasmussen (2010). They put together the Concept Image-Definition along with the emergent model from the RME tradition to create a framework for defining activities. They describe four types of activities in the emergent model: Situational, Referential, General, and Formal. A situational activity is one in which students are familiar with the objects they are studying. A relational activity is one in which students are using something they have learn outside of the area they initially learnt it. A general activity is mathematical work within a newly constituted reality, where interpretations and solutions are within that reality (however, thinking about things can involve other contexts). A formal activity is where the reasoning is done completely within the newly constituted reality. The other dimension of their framework consists of: creating a concept definition, using a concept definition, creating a concept image, using a concept image, and creating a mathematical reality. Assumption digging involves situational activities since students are familiar with the objects they are engaging with. It also involves formal reasoning within a particular reality. However, while students are not creating a completely new reality, they are creating things within that reality which allow them to get a better understanding.

Zazkis & Leikin (2008) worked with teachers on the definition of a square. While this was done with teachers, it would be reasonable to assume that some of the insights are applicable to students. The insight most closely related to Assumption Digging was that some teachers defined squares in terms of other shapes (such as rectangles, rhombuses) while others thought of that as 'cheating.' This is an interesting since they didn't seem to see defining a square in terms of sides or straight lines as 'cheating.' However, from the point of view of the mathematics, these are essentially the same.

To summarize, the examples cited above point to the importance of the relationships between concepts of a theory with other concepts in the theory as well as with the axioms of the theory. It is these relationships which Assumption Digging into a particular theory is looking to explicate. In order to see whether it is successful, it is important to figure out the various ways students engage with Assumption Digging, especially in relating different aspects of the theory.

### *Proof.*

The literature on proof is vast. For the purposes of this review, I will focus on a few of the more influential pieces in this area which are also relevant to Assumption Digging. Harel & Sowder (1998) introduced the concept of a Proof Scheme. A person (or community's) proof scheme consists of what constitutes ascertaining and persuading for that person (or community) (Harel & Sowder, 2007, p. 7). They classify Proof Schemes into three broad types: External Conviction, Empirical, and Deductive. As Weber, Inglis, & Mejia-Ramos (2014) point out, the goal is not to get students to reject the first two proof schemes since even professional mathematicians do not. However, the goal is to get students to see the value of the deductive proof scheme.

Within the deductive proof scheme, Harel & Sowder describe two sub-types: Transformational and Axiomatic. The proof scheme required for successful Assumption Digging is the Axiomatic Proof Scheme.

Tall (Tall et al., 2011), while talking about the cognitive development of proof, talks about three distinct forms of development:

1.  The embodied, which is to do with conceptual understanding

2. The symbolic, which involves translating operations into symbols, and manipulating them

3. The formal, which involves the development of axiomatic formalism, which is in terms of something like set theory

These are not disjoint. For instance, Tall says that Euclidean Definitions and Proofs exist in the intersection of the Embodied and the Formal. Assumption Digging in Euclidean Geometry, hence, is also in this intersection. There is not really any symbolic manipulation involved. We could imagine Assumption Digging exercises in areas like Number Theory, where the symbolic aspects would play a major role.

What is pointed to in the literature is the need for students to be able to reconcile formal and conceptual aspects of proof. While it is true that mathematicians do not always use a deductive proof scheme to convince themselves, it is important that students see the value of a deductive proof. While a deductive proof scheme is required for students to engage in Assumption Digging in a meaningful way, it is worth exploring as to whether Assumption Digging can help develop this appreciation.

**Euclidean Geometry**

Gonzalez and Herbst (2006) identify four main arguments given historically for the existence of the geometry course in the high school curriculum in the US: the formal argument, the mathematical argument, the utilitarian argument, and the intuitive argument. The first two are the ones most closely related to this thesis. The formal argument is that the Euclidean Geometry course helps develop deductive reasoning abilities, while the mathematical argument is that the course allows students to experience what it is like to be a mathematician. While not specifically referring to geometry, the position paper on mathematics in the Indian National Curriculum

Framework points to the need to engage in the discipline of mathematics, and learn logical thinking through mathematics (National Council of Educational Research and Training (India), 2006). If we accept these goals, we are bound to the conclusion that students who go through the course should have a rigorous understanding of the relationships between the definitions, axioms, undefined entities, and theorems of Euclidean Geometry.

However, many have criticized the Geometry course in the US over the years as a 'caricature of authentic mathematics'(Christofferson, 1930; Usiskin, 1980; Weiss & Herbst, 2015; Weiss, Herbst, & Chen, 2009). There are too many postulates, lack of clarity, and a focus on form rather than substance. This is what makes Euclidean Geometry an ideal space for Assumption Digging. If Euclidean Geometry were to be taught in a rigorous manner, with students having a clear understanding of the relationships between axioms, definitions, undefined entities, and theorems, the Assumption Digging exercise would just be students repeating what they have already learnt. That assumption digging turns out to be a valuable exercise could indicate that at least some of this is lacking in the existing Euclidean Geometry course or in its teaching.

**Van Hiele Levels.**

The Van Hiele Levels is an influential model for Geometry learning. There are five levels, usually numbered 0-4 (Vojkuvkova, 2012). The first level is recognition of geometrical objects. The second level is to do with describing shapes. The third is being able to abstract away from specific shapes to categories and properties. The fourth involves an understanding of the distinction between axiom, definition, theorem and proof, and students are able to differentiate between necessary and sufficient conditions. The final level is an understanding of how mathematical systems are established, and students are able to do all types of proofs. They are also able to comprehend Euclidean and non-Euclidean geometries and their relationships to the axioms. This Assumption Digging exercise resides somewhere in the 4th and 5th levels. Unearthing the relationships between different statements of the theory is one of the goals of Assumption Digging. Whether Assumption Digging helps students to internalize these relationships, both at the level of Euclidean Geometry and at the level of the nature of a Mathematical Theory, would be a valuable exploration.

**Research Questions**

There are many interesting avenues for research into Assumption Digging suggested by the literature: the relationship between Assumption Digging in Euclidean Geometry and the Van Hiele levels, whether Assumption Digging helps develop a deductive proof scheme in students, the relationship between Assumption Digging and the reconciliation between the conceptual and figural aspects of Geometry, amongst others. Each of these would constitute a research project on their own. I will address these issues as and when they come up in the data. However, in order to even begin thinking about these things in a clear and rigorous manner, we need to explore some far more basic questions about the nature of interactions in a classroom session on Assumption Digging. The following are the research questions I'm intending to investigate:

1. What are the ways in which students engage in Assumption Digging in Euclidean Geometry?

2. What are the ways in which the instructors' interventions direct Assumption Digging sessions?

# METHODOLOGY

## Instructional Setting and Student Background

### The Workshops.

The relevant sessions were carried out during workshops in two schools in Pune, Maharashtra, India over the summer of 2018. The workshops were on mathematical thinking more generally, but with a focus on Theory Building. I will be referring to the two schools as Indus School and Ganga School. They are low cost private schools, owned by the same non-profit.

The workshops were with students from Grades 7-10, with the median student being from the $9^{th}$ grade. The main goal of the workshops was the professional development of the teachers. The school administration's goal was that the teachers move from developing mathematical understanding in students to developing their ability to think like mathematicians. The sessions were recorded for teachers who were not able to attend due to other commitments.

The workshops for each of the schools were four days long, and the sessions on Assumption Digging happened towards the end of the workshops. The lesson plans for some of the other sessions are in the appendices. The following is the list of the workshop sessions.

Day 1 (Sat – Indus): Main activities – Straight Lines and Intersections, Circumscription (Assumption digging into Definitions)

Day 2 (Sun – Indus): Main activity – Discrete Geometry (Theory Construction)

Day 3 (Mon – Ganga): Main activities – Straight Lines and Intersections, Circumscription

Day 4 (Tue – Ganga): Main activity – Discrete Geometry (Theory Construction)

Day 5 (Wed – Ganga): Main activity – Coloring Problems

Day 6 (Thurs – Ganga): Main activity – Assumption Digging in Euclidean Geometry

Day 7 (Sat – Indus): Main activity – Assumption Digging in Euclidean Geometry

Day 8 (Sun – Indus): Main activity – Coloring Problems

**The Facilitators.**

The workshop was mainly facilitated by two people – a retired theoretical linguist (referred to as Raghavan) and me, a mathematics graduate student. Apart from that, there were two others who helped, who will be referred to as Mira and Rohit. Mira is an educator who teaches occasionally at Indus. Her main interactions are with 12th grade students. Rohit was 21 years old at the time of these sessions. At the time, he was in the process of applying to colleges for his undergraduate education. He was homeschooled and has been taking help from Raghavan and me with his education. Mira and Rohit's main interventions were related to classroom management. However, Rohit did occasionally intervene in relation to the actual content of the session.

**Education in India**

In India, full-time government or private school education is mandatory for students from ages 6 to 14. Schooling children solely at home is also permissible, but the vast majority of students of all economic classes attend either government or private school. Government schools are free of cost, and there are both low and high cost private schools. Fees at private schools are typically paid by the family of the student. However, in some parts of the country, the government pays for some economically disadvantaged students to attend private schools. Some private schools have religious affiliations (such as Hindu, Islamic, and Christian schools). The language of instruction varies from school to school. Some schools teach in English, others in Hindi, and still others in various local languages.

The school year is different in different parts of the country. In Maharashtra, it starts in June and ends in April. Students from upper and middle class background typically take various classes over the summer holidays.

Schools are usually affiliated to an educational board. There are two national boards in India – Central Board of Secondary Education (CBSE) and Indian Certificate of Secondary Education (ICSE). Apart from that, most of the states have their own boards. Associated with the boards is a textbook writing and teacher education institution. For the CBSE board, that role is played by The National Council of Education Research and Training (NCERT). While there is a prescribed syllabus for each subject for each grade, schools have the freedom to do what they want on a day-to-day basis as long as they prepare their students for the board exams in the $10^{th}$ and $12^{th}$ grades. There are also schools which are affiliated with the International Baccalaureate (IB) and the International General Certificate for Secondary Education (IGCSE). These schools tend to be much more expensive.

At the elementary school level (ages 4-10), students are usually taught all the academic subjects by a single teacher. From middle school onwards, teaching is more specialized. Most teachers are required to have a Bachelors in Education (B Ed), with high school teachers requiring a Master's Degree in their discipline. This is granted by various universities, and the curriculum varies across them.

Students do not get a choice in what subjects they study till the end of the $10^{th}$ grade, apart from some choice in learning languages. Every student has to take the same mathematics course till the $10^{th}$ grade. After that, they can choose to drop mathematics. The mathematics till the $10^{th}$ grade includes arithmetic, algebra, geometry (Euclidean and analytic), trigonometry, linear algebra, and probability. For those who take mathematics till the 12th grade, the main above topics are covered in more detail, with the addition of calculus.

The 10th and 12th grade examinations are run by the respective boards. They usually happen outside of the school premises. The 10th grade CBSE examination had been stopped by the previous central government in favor of continuous assessment. It was reinstated for the 2017-18 academic year. The questions in the exam do not allow for much time to think due to the number of questions. The focus, especially from the 9th to the 12th grade in most schools is on doing well in the exams. Hence, understanding is ignored in favor of quick strategies to solve test problems.

**Student Background.**

The schools are mid-priced private schools affiliated to the CBSE board for the purposes of 10th and 12th grade examinations. For mathematics, they use textbooks from NCERT, and the medium of instruction is English. While the textbooks can possibly be used to create a classroom environment as suggested by Engle and Conant, the textbooks themselves don't encourage the four principles mentioned in the Introduction to this thesis. For instance in the Euclidean Geometry chapters, students are given definitions, axioms, and proofs to all the important theorems. They have to come up with proofs themselves, but only to less crucial theorems. Conjecturing is not explicitly encouraged.

In the session in Indus, there were 12 students, 5 girls and 7 boys. In Ganga, there were 23 students, 12 girls and 11 boys. The students are mostly from Middle Class families from Maharashtra. A few of them have had similar previous mathematical experiences to this session – not assumption digging, but things such as defining, reasoning, etc. This happened during previous workshops with Raghavan and Mira. In terms of geometric knowledge, all of them have previously had an introduction to straight edge and compass constructions, and most of them, those in the 9th and 10th grades, have had some introduction to Euclidean Geometry.

**Instructional Setting.**

*Indus.*

The session was carried out in a large room with two large green boards, giving substantial room to leave bits up for the entire session and refer back to them. The students each had notebooks and writing implements. However, none of them had access to computing devices. There were two types of pedagogies used – whole class discussion and group work in groups of 3-5. Students would typically break into groups to discuss things where they were encouraged to share ideas, and then contribute to the whole-class discussion. While students ranged from grades 7 to 10, I only know the specific grade levels of two of the students – Tushar was in grade 7 while Tarini was in grade 9.

The workshop sessions were during the weekend, and students volunteered to take the workshops. The students who were in grades 7 and 8 had previously taken a year-long course in their 6th grade for which the material was provided by Raghavan and me, amongst others. This course was on inquiry and critical thinking. It involved mathematics, but was not restricted to it. A few of the older students had been part of one-off workshops on inquiry and critical thinking, which Raghavan had run in the school. Some of the other students had heard of us and had probably seen Raghavan around due to his association with the school.

*Ganga.*

The room used in this school was small and crowded. There was almost no space to move around, especially for instructors to get to the back of the room. There was also only one small green-board. Hence, it was not possible to have statements on the board for the entire session, and referring back required re-writing or was done verbally. There was also whole class discussion and group work, but the lack of space restricted the groups and especially interactions between groups. The sessions here were conducted during school hours, which resulted in teachers being present for small portions of the workshop, and not participating. While the students ranged from grades 7 to 10, I do not know the specific grade levels of any of the students. Only a few students in Ganga had attended one-off workshops run by Raghavan; however, some others knew of us by reputation.

**Prior Geometry Experience.**

Euclidean geometry is introduced more formally in the 9th grade in the NCERT textbooks. Prior to that, the geometry chapters in the 7th and 8th grades are concerned with compass and straight edge constructions, and with understanding polygons and their properties.

The Euclidean Geometry course begins with an introduction to Euclid's axioms, postulates, and definitions. At the end of that section, it states that while what comes next is Euclidean Geometry, it isn't based on the same axiomatic system as Euclid's due to it not being satisfactory.

**The Mathematics.**

It is important to get a grasp of the mathematics itself. The claim made in both the sessions was: Any two triangles with the same base and between the same parallel lines have the same area. In both the sessions, students provided the initial proof for this claim – they had seen the proof before.

Here is the proof:

Premise 1: The height of two triangles between two parallel lines is the same

Premise 2: The area of a triangle is ½ base x height

Premise 3: The bases of the triangles are the same

Conclusion: The areas of the triangles are the same

Premise 3 is given in the claim, so the only premises to investigate with Premises 1 and 2. Both sessions investigated both the premises.

## Data Collection

The data consisted of the video recording of the sessions – around 45 minutes in Ganga and 1 hour in Indus. The intention of the video recordings was for the professional development of the teachers. The teachers who weren't present would have the ability to go through the videos, and those who were present could use them as a reminder of what went on in the sessions. Hence, given monetary constraints, it was determined that a one-camera set up focused on the green-board would be sufficient. So, students are only visible when they come up to the green board. Also, since the mic was aimed at the green board, parts of the students' speech weren't completely clear.

Portions of the videos of the sessions were then tagged with which aspects of Assumption Digging students engaged with. The tagging was in the form of AD1, AD2, AD3, and AD4, corresponding to the aspect of the definition of Assumption digging, as shown below:

AD1:   Given a claim, extract possible premises on which the claim relies.

AD2:   Attempt to clarify the meanings of the words used in the claim and premises.

AD3:   Continue 1 and 2 recursively

AD4:   At some point in this recursive exercise, state the premises as axioms and definitions, and decide that a particular word is an undefined entity with certain properties (which also form some of the axioms of the theory).

A particular portion could have more than one tag. These portions were then transcribed. This transcription gave the gist of what went on in the sessions and was not an exact replication.

**Data Analysis Method**

The data will be presented in a chronological manner for each of the two sessions. The reason for doing this is that it will give the reader a sense of the direction of the sessions. Especially while talking about the second session, I will attempt to point out the ways in which it differed from the first. I picked those portions out of the previously transcribed ones such that there was something which distinguished that from the previously chosen portions (in terms of the two foci of the analysis mentioned below). These were then transcribed in detail.

Students were given pseudonyms. Since the students are not always visible in the footage, there might be more than one pseudonym for a given student.

Given that the goal of the analysis is to address the research questions, there are two main foci. The first is the various ways students engaged in Assumption Digging. While the tagging to the aspects of the definition sheds some light on this, the goal of the analysis was to uncover the ways they went about engaging with the different aspects of Assumption Digging. While doing this, I was looking for instances which revealed students' proof schemes (Harel & Sowder, 2007) and their relationship to the Assumption Digging. I was also looking for the disconnect between the 'conceptual' and 'figural' aspects of Geometry (Mariotti & Fischbein, 1997) where and if applicable and see how they were addressed, either by Assumption Digging or in order to move forward with Assumption Digging. In a few situations, I have also attempted to represent student arguments in the form of tree diagrams. Since the tree diagrams are less useful when the flaw in the reasoning is in the form of a hanging leaf not stated as an axiom or a definition, I have done this when other types of flaws emerged.

The other focus is the role of the instructors. As in the RME tradition, the role of the instructor is expected to be that of a facilitator. Hence, while analyzing the data, I looked for the various ways in which the instructors intervened, how that shaped the discussion, and pointed out various alternative routes they could have taken students down. I was also looking at the extent to which the discussions were shaped by the students, and the extent to which they were shaped by the instructors.

# RESULTS

## Preliminaries

It would be useful to take another look at the definition of Assumption Digging I will be using.

AD1:   Given a claim, extract possible premises on which the claim relies.

AD2:   Attempt to clarify the meanings of the words used in the claim.

AD3:   Continue 1 and 2 recursively

AD4:   At some point in this recursive exercise, state the premises as axioms and decide that a particular word is an undefined primitive with certain properties (which also form some of the axioms of the theory).

The analysis will be in two sections. The first section is about the session of Assumption Digging in Ganga, and the second is about the same session in a school in Indus. The analysis below is chronological in order to get a sense of the direction of the sessions. It also highlights the multidimensional nature of the exercise, with plenty of directions available at any given point. Especially in the second part, I will be highlighting the differences between the two sessions.

## Session 1 – Ganga

In Ganga, the session started with Assumption Digging into premise 2, asking why we should believe that the height of triangles between two parallel lines remains the same. The following is an exchange inquiring into what we mean by 'height of a triangle.' MK refers to me and Raghavan is my fellow instructor.

### Ganga: Altitude of Triangle (AD2).

*Arean: And then. height is.. height also remains throughout. In all the //triangles.*

*Raghavan: Height//*

*Raghavan: No. For example, take this one. This is that height and this one is this height. It's different.(Board not visible)*

*Arean: No. Uh. Height is taken from the uh vertex.*

*Raghavan: This is not the height? //Where is it?*

*Ananya: The altitude//*

*Raghavan: That's the altitude, right?*

*Students (multiple): No*

*Ananya: Starts from the vertex*

*(Students responded together, mostly inaudibly. Ananya was the only one audible)*

*Raghavan: Starts here. Okay. Wait, you're saying this is the same as that?*

*Students (multiple): Yes.*

*MK: Why?*

*Raghavan: I don't believe that.*

When Arean claimed that the height of all the triangles was the same, Raghavan picked a random spot on one of the sides of the triangle and dropped a perpendicular. Clearly, the lengths of these was not the same across the triangles. What was being pointed at was the need to clarify the definition of 'height of a triangle,' and bridge the disconnect between the conceptual and figural aspect of the concept. This is the first example of AD2, and gives evidence that assumption digging of this type is possible for students to engage with. Immediately Arean responded saying that height is from the vertex and not from a random point on the side, and Ananya introduced the concept of altitude. Apart from giving the initial proof, which could be argued was a re-telling of what they had previously learnt, this was the first example of Assumption Digging in the session.

The discussion continued with the students attempting to address the skepticism of the instructors:

**Ganga: Parallel Lines (AD1 & AD2).**

*(Straight after the previous extract)*

*Students (multiple): Because they are parallel lines.*

*...*

*MK: What is a parallel.. What are parallel lines?*

*(Laughter from students)*

*Gayatri: The lines who do not intersect each other and have a equidistance between...*

*Raghavan: Okay, lets write down the definition*

*(Students talking in the background inaudibly)*

*Gayatri: ... on the same plane.*

*Raghavan: So we had three, three definitions of parallel lines here, and you are choosing this one, right?*

*Students (Multiple): Yes.*

*Raghavan: And you are choosing this one. Parallel lines, when extended, do not intersect.*

In the first part of this exchange, some students attempt to justify their previous claim by saying it is based on the premise that the lines are parallel, indicating that this is an example of AD1. The clarification of the definition of parallel lines is an example of AD2.

Raghavan questioned their claim that all the altitudes have the same height, and they attempted to justify it by appealing to the lines being parallel. They are making use of a deductive proof scheme here for the few steps of reasoning in the above exchange. However, they didn't explicitly share their steps of reasoning from parallel-ness to the conclusion.

By asking 'What are parallel lines?', I am pushing them towards that. An alternative way for the instructor to ask the question would be: 'Why does the lines being parallel result in the conclusion that the altitudes are the same?' That would have required the students to figure out for themselves that they need to state what they mean by parallel lines. At the end Raghavan mentions the 'three definitions' of parallel lines. This is appealing to a previous session with the students during the workshop when they defined parallel lines in three ways (the laughter probably was as a result of the memory of the previous session. The lesson plan for that session is in Appendix 1). These definitions were:

1. Straight lines which are equidistant (students re-state the definition of equidistant in the next section of dialogue)

2. Straight lines which do not meet

3. Straight lines A and B are parallel if there exists a straight line C which is perpendicular to both

**Ganga: Equidistant (AD1 & AD2).**

*MK: What… What was the definition of equidistant?*

*(Referring back to previous session)*

*Joseph: When 2 points have same distance throughout?*

*MK: What do you mean 2 points? Two points always have the same distance.*

*...*

*Jaya: Whether you take more.. like if we take 10 cm on the segment, the distance will be the same (followed by something inaudible)*

*MK: What do you mean the distance? From..*

*Jaya: (Something inaudible)*

*MK: Distance from start to end point?*

*Jaya: Na the point on CD and AB (See Figure 6)*

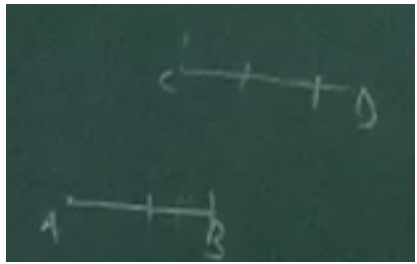

*Figure 6.* Screenshot of green board from conversation with Jaya.

*MK: (Something inaudible while handing the chalk to Jaya)*

*(Jaya proceeds to the board to explain. Her explanation is mainly in gestures, but MK repeats it in words after)*

*MK: And you move along the same distance here and here. And this will be the same as… as this. Okay.*

*(Raghavan's phone rings while MK speaking, and he leaves the room)*

*MK: Okay. So, now, how do we know that is happening here? Your.. so, say we start here, and we have this distance from here to here and we move from here to this point. How do we know that this distance is the same as this distance?*

*Abdul: Because.. um.. the.. um.. lines which are parallel… Can I draw it and say.*

*MK: Yeah.. Let's.. Let me remove the.. for now we only care about the.. We only care about this.*

*(Raghavan returns to the room)*

*(MK Rubs the triangles on the parallel lines and draws only the two perpendiculars – see Figure 7)*

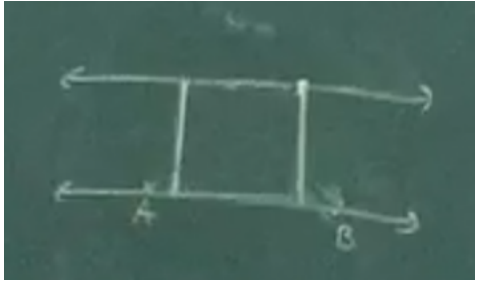

*Figure 7.* Screenshot of green board from conversation with Abdul.

*Abdul: Um.. These 2 lines are parallel, and these 2 lines are parallel, and um.. these 2 points are um.. equidistant.*

*MK: Ahh.. So, you are saying that these two lines are also parallel. How do you know that?*

*Abdul: Because they form a ninety degree angle.*

*MK: Where is that in equidistant? (Pointing to the definition of parallel lines as equidistant)*

*(Abdul points to the third definition on the board)*

*Students (multiple): Common perpendicular (the third definition)*

At the beginning of the discussion, I ask what the definition of equidistance is. The students' engagement with that is an example of AD2. Then, they use the definition as a premise in the justification of another claim, which is an example of AD1.

Rather than directly asking what equidistance is, another way for me to convey something similar would be: 'How does equidistance of parallel lines give us that the altitudes of the triangles between a pair of parallel lines are the same length?' This may have resulted in the students realizing the need for the definition of equidistance (which they had already come up with in a previous session).

The students involved in this discussion used two of the definitions of parallel lines. During a previous session, they had seen that these definitions were equivalent on the plane. By the end of the discussion, they have given a reasonable justification that the altitudes of the triangles have the same length. Jaya draws a perpendicular connecting two of the parallel lines. She says that what makes the two lines equidistant (and hence parallel) is that when you move the same distance along both the lines in the same direction starting at the perpendicular, and connect the end points, you get a line of the same length as the perpendicular. This new line is also perpendicular to the two parallel lines, which Abdul suggested was due to the third definition. This argument, which makes use of the equivalence of two definitions of parallel lines, is quite sophisticated and gives more evidence that at least some of the students are making use of a deductive proof scheme.

The class then moved on to Premise 2. They initially started by attempting to justify the area of a triangle using the areas of rectangles and general quadrilaterals. However, they quickly moved on to parallelograms.

**Ganga: Breaking a Parallelogram into Triangles (AD1).**

*Aditya: Suppose if we take a parallelogram, the area of parallelogram is base into height. So, there are two parts there and can be two triangles. If we divide it by a diagonal, there are two triangles. So it can be half base into height.*

*MK: So, you are say that the triangle… Ok, so there is something called a parallelogram. (While writing on the board)You are starting with the assumption that area of.. of parallelogram is base multiplied by height. You can break parallelogram into 2...*

*Aditya: Into two triangles.*

*MK: Ya.. How do you do that?*

*Students (multiple): Diagonal*

*MK: With a diagonal.. Okay.. Then?.. Then?*

*Aditya: Then the area will be half base into height. Because total parallelogram is base into height.*

Even when justifying the area of triangles with quadrilaterals and rectangles, they were moving in a productive direction – attempting to justify the area of a triangle in terms of areas of other shapes. This was their first instinct when presented with the question. This is very different from the experience Zaskis & Leikin (2008) had with teachers, many of whom thought of this strategy as 'cheating.'

Aditya's argument is attempting to justify the claim that the area of any triangle created by cutting a parallelogram into two via the diagonal is ½ base x height. He does this by appealing to the area of the parallelogram and the result of a construction, which is an example of AD1. However, just because the triangles you get by breaking parallelograms have the property we want, that doesn't mean that given any triangle, it has that property. It appears that they didn't automatically see this distinction, and it had to be explicitly stated to them. In the following conversation, a student (Bharat) attempts to fix this after I point out the problem.

### Ganga: Making Parallelograms from Triangles (AD1).

*MK: Does someone want to explain the question?*

*Bharat: If you draw random triangles how do you get area is half base into height?*

*MK: This is a proof given on.. for triangles which are formed by cutting parallelograms into two. Uh yeah..*

*Bharat: If you take the.. a triangle, and turn it.. the same triangle, and join the bases, you get a parallelogram*

*MK: So, you're saying every triangle..*

*(MK erases board, draws a triangle)*

*Raghavan: Make a parallelogram*

*Student draws on the board.*

*MK: What do you need to do to do this?*

*Raghavan: How do you create a parallelogram from the triangle?*

*MK: What is the procedure?*

*(Student makes some gestures which look like flipping the triangle across a side)*

*Raghavan: Flip it*

*Bharat: Flip it (See Figure 8)*

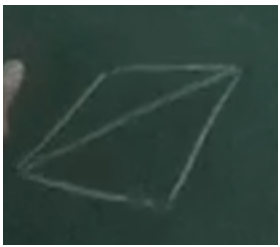

*Figure 8.* Screenshot of green board from conversation with Bharat

The session then moved to showing that given a triangle, it can be made into a parallelogram. Bharat was taking the first step towards extracting premises in order to justify that claim, an example of AD1. The initial suggestion by Bharat works perfectly. Rotating a copy of the triangle and sticking it to the triangle results in a parallelogram. However, when showing it on the board, he appeared to indicate, via a gesture, that the triangle was to be flipped. Once Raghavan suggested flipping, Bharat accepted. Even though the flipping would actually work on the example triangle drawn by me (which was isosceles), it doesn't work generally. There was a clear distinction between the conceptual and figural aspects, which was brought about by considering a very specific example. However, both Raghavan and I accepted it and actually thought that it worked generally. Raghavan and I being wrong resulted in the following valuable discussion.

**Ganga: Student finds Flaw in reasoning (AD1 & AD3).**

*Manu: He said that base of the triangle can be any side. And as we see as base.. it becomes a triangle.*

*MK: Can you show me what you're saying.*

*(Manu draws a larger triangle from a triangle by flipping on the board and gives an explanation of what he is doing, which is not completely audible – see Figure 9)*

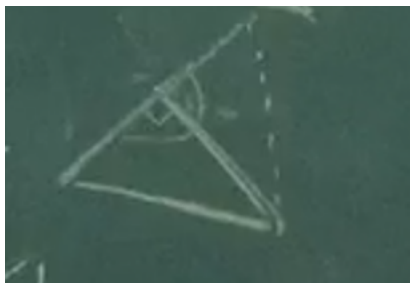

*Figure 9.* Screenshot of green board from conversation with Manu

*MK: Okay. So, he is saying you take this and attach it like this, it becomes another triangle.*

*Manu: Yes.. Base can be any side.*

*MK: Ya. Base can be on any side.*

*MK: So, this becomes a triangle*

*Manu: Yes*

*MK: Why? Why would it be a triangle?*

*Manu points at its shape.*

*MK: No. No. That it looks like a triangle is not good enough. Why is it's a triangle?*

*Ok, let me draw it for you. What about this triangle make it that when you flip it, it is not a.. it is still a triangle?*

*Manu: (Something inaudible)*

*MK: No.. No.. That yeah. So it remains three sides. But what about this tri.. So, if I take this triangle and flip it, you have a quadrilateral. It looks something like this. But this triangle. Why is not it? Why is this a triangle and that not? (silence for a few seconds) Ok, why does it remain a triangle… Anyone.. wants to help him?*

*MK: (pointing at a student with their hand raised) Okay. Suggest something to him..*

*Vidur: Only right angled triangles..*

*MK: If this is a right angled triangle and you flip it, why does it remain a triangle?*

*Vidur: Because this is straight. If we flip a congruent triangle, it will still be straight.*

*MK: Why?*

*Students (multiple): It's 90*

*Vidur: Because these two triangles are congruent.*

*MK: Ah.. Because these two triangles are congruent, this becomes two right angles, and as we talked about that day, is a straight angle, right? You.. you.. you guys told us what a straight angle is.*

*Raghavan: Two right angles. That's a straight line.*

*MK: So, that means this will be a straight line.*

*...*

*MK: He is pointing out a flaw in this argument.*

*Raghavan: Ah.. I see*

*MK: That, when you flip it, you don't necessarily get a quadri.. a parallelogram.*

*Raghavan: You might get a triangle.*

*MK: Yeah.*

*Raghavan: Ah.. Okay. The argument was you can create a parallelogram by flipping it. It will always be a parallelogram. He is saying that is not true. Okay.*

*MK: Not by flipping.*

*Students (multiple): We can still make a parallelogram.*

*MK: No..No..No.. So, is flipping the right procedure?*

*Students (multiple): No.*

*MK: Then, what would it be?*

*...*

*Serena: So, in the argument, he said you should flip it, but he didn't say that you should flip it from the longest side so that it will form a.. a parallelogram.*

*MK: Okay.. Let me.. let me try giving you another triangle.*

*Raghavan: Did everybody get her.*

*MK: So, when we flip it, you get this, right? (See Figure 10)*

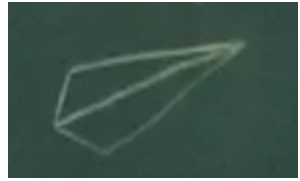

*Figure 10.* Screenshot of green board from conversation with Serena

*Serena: We get quadrilateral*

*MK: Yeah. We get a quadrilateral. But, we don't get a parallelogram.*

Figure 11 represents an explicit version of the argument Bharat made in the previous exchange. What Manu is doing here is questioning the steps of reasoning which take you from the premises to the conclusion. He is implicitly questioning the existence of Premise 1.

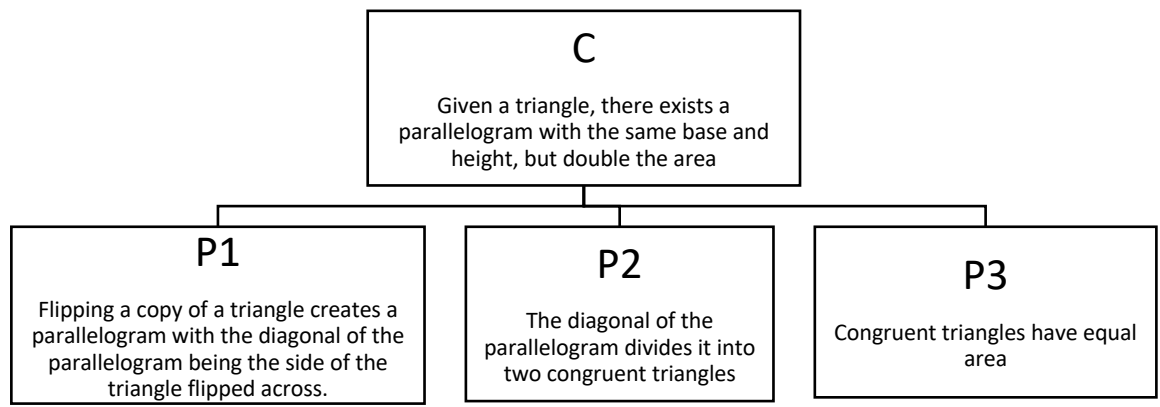

*Figure 11.* Tree representation of Bharat's argument

If you take a right angled triangle and flip it across one of its non-hypotenuse sides, and glue it back to the original, you get another triangle rather than a parallelogram. This is what Manu was arguing. In this case, without prompting, Manu found a flaw in the reasoning the class had agreed on, including the instructors. He was able to point to the need to reconcile the conceptual and figural. It is not the case that by flipping a triangle across a side, you always get a parallelogram. Serena then goes on to try and fix the problem by suggesting that you can flip the triangle across another side. This is an interesting suggestion since it works for some classes of triangles such as isosceles right angled triangles. However, for other triangles, you get a general quadrilateral. As Serena correctly pointed out, this quadrilateral need not be a parallelogram. This was the first time during this session that any student saw a fault in the reasoning without being prompted. This exchange gives more evidence that students are using a deductive proof scheme since they are willing to accept a claim being wrong when presented with just one counter-example.

In this exchange, Manu was asking why we should believe a particular claim (the P1 from the diagram above), an example of AD3. In doing so, he found a counter-example which rejects P1, which is an example of AD1. You could also think of this as questioning C. In that case, Manu has detected a flaw in the edges of the graph, the steps of reasoning. With that perspective, the edge associated with P1, and hence P1, needs to be removed and replaced with some other premise. This exchange shows that students are able to initiate assumption digging, at least in situations where there are obvious flaws in the edges.

The exchange with Vidur was tangential to the rest of the discussion, but led to an interesting theorem that if you flip a right angled triangle along a non-hypotenuse side, and glue the two triangles together, you get another triangles. The effect of transformations and gluing on shapes would be an interesting session on its own. The hole in the argument was then plugged by introducing a new procedure for converting a given triangle into a parallelogram as shown below.

**Ganga: Fixing the argument (AD1 & AD2).**

*Jaiveer: If you take a triangle, the smaller side. Which will be the smaller side, you will draw parallels line from the.. the opposite vertex.*

*MK: Oh.. Okay.. So, you want to draw a parallel line through this?*

*Jaiveer: Yes and connect (followed by something inaudible)*

*MK: Okay. I'll draw a parallel line and connect it.. (See Figure 12)*

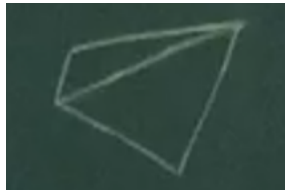

*Figure 12.* Screenshot of green board from conversation with Jaiveer.

*Students (multiple): No, of the same measure. The same measure.*

*MK: Oh. Okay. So, this should be the same measure as this and parallel.*

*…*

*MK: Why is this a parallelogram?*

*Raghavan: A parallelogram should be that these two lines and parallel and these two are parallel. You drew it such that they are parallel.*

*MK: You said this is parallel to this and this length is the same as this. Why does it mean that.. uh.. these two are parallel to each other? Right, for a parallelogram you need that these two are parallel.*

*…*

*Kanika: Because if they are parallel and measure is same, if you connect it, it will be same. Because if it were parallel and measure would be different, we would have got a different (followed by something inaudible)*

*MK: One second.. Use what parallel means.*

*Kanika: Parallel means they are equidistant everywhere.*

*MK: Okay.. So, how does this.. How are this.. How are these 2 lines equidistant.*

*Kanika: These two lines are of equal measure*

*MK: We don't know that these two lines are equidistant. We know these 2 lines are equidistant*

*and of equal length, right?*

*Kanika: Yeah.*

*MK: That's how he made 'em.*

*Kanika: (who has walked up to the board) If this is the distance between these 2 points, it will*

*also be the same distance for these 2 points. (See Figure 13)*

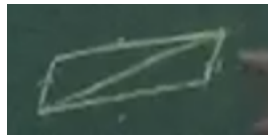

*Figure 13.* Screenshot of green board from conversation with Kanika.

*MK: Ah.. Okay.*

*Raghavan: Good*

*MK: S.. So, since these two are parallel, they are equidistant. So, that means this length is the*

*same as this length.*

*Raghavan: Because equidistant is defined as if you draw a little bit distance here and the same*

*distance here, this will also be the same. That's the way they defined it.*

These students then attempted to fix the problem by using another procedure for converting a given triangle into a parallelogram. The suggestion by Jaiveer was to draw a line parallel to one of the sides of the triangle. Mark off a point such that the length is the same as that of the line it was drawn parallel to, and connect the marked point to that line. His suggestion of using the 'smaller side' was not explored, but could have been. The conclusion relied upon explicating the definition of a parallelogram, an example of AD2. The definition they initially had only said that opposite sides are parallel. Hence, they were pushed to show that their procedure resulted in opposite sides being parallel. In order to do this, Kanika appealed to the equidistance definition of parallel lines, explaining that equal length of two opposite sides means that the other two sides are parallel. She was engaging in AD1.

The rest of the session involved justifying the area of parallelogram in terms of area of rectangle, and at the end the instructors hinted that they could try justifying that in terms of area of a square. The session ended with the instructors pointing to the fact that they could have started with any of the definitions of area (triangle, parallelogram, rectangle), and deduced the rest from that definition.

**Session 2 – Indus**

In Indus, we began with the Area of a Triangle rather than equivalence of heights as we had done in Ganga. The height part was similar to that in Ganga, so I will not be presenting data from there.

**Indus: Area of Triangle (AD1).**

*Tarini: This is a rectangle.. Length and breadth. Area of this will be length into breadth..*

*(Something inaudible). So, if we cut it in half, you get half base into height.. (See Figure 14)*

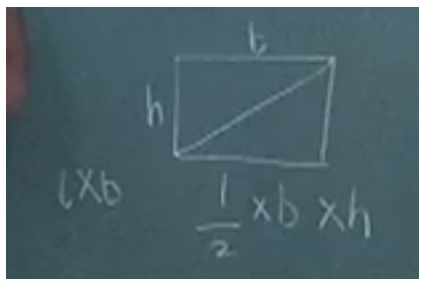

*Figure 14.* Screenshot of green board from conversation with Tarini.

*MK: Okay.. let me write that out.*

*(Rohit says something inaudible in the background)*

*MK: Area of rectangle is base into height. If you break a.. a rectangle in half using a diagonal, the area of each resulting triangle is half the area of the rectangle. So, area of any of these triangles is half base into height.*

After discussing in groups, the first argument a student gives is that the area of a rectangle is base x height, and you get a triangle when you break a rectangle into two using a diagonal. Hence, the area of each of those triangles is ½ base x height. Tarini is extracting premises in order to justify her conclusion, an example of AD1. The problem here is the same as in Ganga – while breaking a rectangle into triangles results in two triangles with the properties we seek, that doesn't necessarily mean that given a triangle, that triangle's area is ½ base x height. Straight after Tarini presented her argument, another student suggested using a square rather than a rectangle. This resulted in the following invaluable exchange:

**Indus: Why not a square? (AD1).**

*A student suggests using a square rather than a rectangle. Most of what that student says in inaudible. Tarini responds.*

*Tarini: So, we did not take a square.*

*Raghavan: Did everybody get the question?*

*Students (multiple): Yes.*

*Tarini: Because, the sides of a square are equal and if we cut that into half, the two sides of the triangle will be equal and one side will not be.*

*The same student from earlier ask again about the square. Most of what he is saying is inaudible.*

*Tarini: No.. They want a generic triangle*

This is an extremely interesting exchange. Tarini sees that what we are looking for is a justification that works for what she calls 'generic' triangles, and not for the particular type of triangle created from a square. She points out that by breaking a square via its diagonal, you only get triangles with two sides equal and not a 'generic' triangle. By 'generic,' from what she says earlier, she appears to mean triangles which you get from cutting a rectangle in half. One possible interpretation is that she doesn't see triangles which cannot be created this way as 'generic.' Another is that she believes that all triangles can be created in this manner.

Rather than continuing on whether any triangle can be made via this procedure, the session took a slight detour to explore the notion of congruence.

### Indus: Congruence (AD2).

*Raghavan: There is one question I would like to ask. I don't believe this. I don't believe if you draw a diagonal*

*MK: it breaks it onto two..*

*Raghavan: That it breaks it into two equal halves. I know that you are surprised. But, what mathematics is about is questioning the obvious. This might be obvious to you, but you have to prove it.*

*MK: You have to show.. you have to show that they have equal area.*

*Avir: that every congruent figure.. if they are congruent, they have same area.*

*MK: What do you mean by a congruent figure?*

*Raghavan: Same dimension? Do you remember the definition of congruence?*

*Avir: If all the angles and sides are equal.*

*MK: So, okay. So, you can't say two circles are equal.. are congruent?*

*Raghavan: If I draw a figure like this, and I draw another figure here, how would you decide whether those two figures are congruent?*

*(Multiple students say something)*

*Raghavan: Put them on top of the other and if they perfectly coincide, they are congruent.*

This was a discussion which was much shorter in Aundh, and the instructors didn't really push the students to clarify the concept of congruence, an example of AD2. Even here, there was no justification given for the claim that the two triangles had equal area. However, the notion of congruence was at least touched upon. The students who spoke seem to have a figural notion of congruence. However, their conceptual notion doesn't seem well developed. This could be an interesting path to go down in future iterations of this session. Avir's claim that 'if two shapes are congruent, they have the same area,' could be taken as the basis for a definition of a general notion of area. Straight after this, the session returned to addressing the procedure Tarini had set out.

**Indus: Do rectangles work? (AD1).**

*MK: Are you saying that any triangle can be made like this. By breaking some rectangle into half?*

*Tarini: If you have a triangle. Then, take the same triangle and rotate it. Then we get a rectangle.*

*MK: Does everyone agree?*

*Students (multiple): No.*

*Then, given a non right angled triangle by MK, Tarini comes to the board to attempt to make it into a rectangle. What she draws is a non-rectangle parallelogram and not a rectangle, and she realizes she was wrong.*

In the previous exchange with Tarini about squares, she had rejected using squares since that would result in non-'generic' triangles. As previously stated, one interpretation of her use of 'generic triangle' is that she doesn't see triangles which cannot be created by cutting rectangles in half as 'generic.' Another interpretation is that she believes that all triangles can be created in this manner. This exchange gives us some reason to believe that the second interpretation is more likely than the first since she did accept a fault with her argument when she was not able to create a rectangle from the given triangle.

After this, another student (Tushar), attempted to clarify the entire argument made so far, by laying out the basic assumptions.

**Indus: Starting from Area of a Square (AD4).**

*Tushar read out an argument. With slight clarificatory questions from me, here is his argument*

*as written on the board by me:*

*Area of a square = 1*

*A triangle has 3 sides and 3 non-straight angles.*

*Any rect can be broken into b x h unit squares*

*So, area of rect is b x h*

*Rect divided by its diagonal gives 2 cong tri.*

*So, area of those is ½ bh*

Figure 15 is a representation of Tushar's argument.

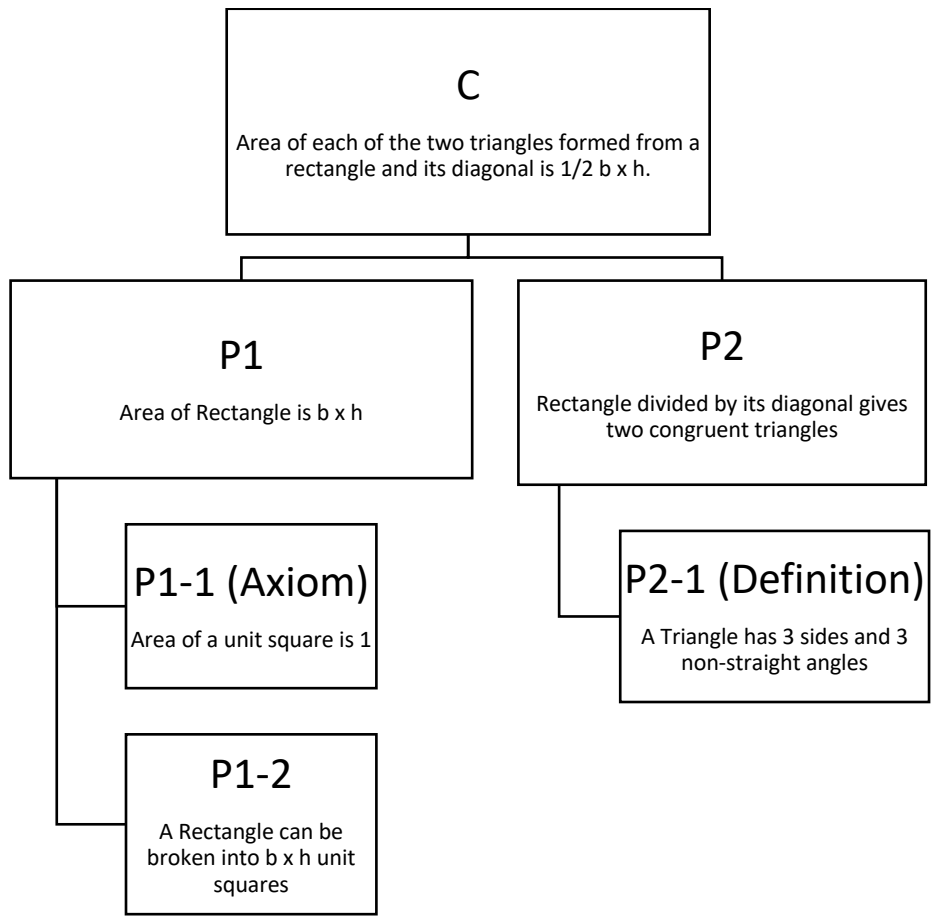

*Figure 15.* Tree representation of Tushar's argument

The interesting part of this is that this is the first time any student laid out an argument explicitly stating at least some of their most basic assumptions. It shows that Tushar sees some value in an axiomatic deductive proof scheme. Although the language was not used, Tushar was setting up the building blocks of the theory – the definitions and axioms. It could have been interesting for the instructors to either spend more time on this or return to this at the end of the session and point out explicitly what Tushar was doing. There is a lot missing in this argument. For instance, P2 doesn't follow from just P2-1. While the designation of P1-1 and P2-1 as axioms and definitions respectively seems reasonable given his argument, it isn't clear as to what designation P1-2 ought to have. However, this is the only example of AD4 in the two sessions. No matter what the flaws, it shows that it is possible for students to engage with assumption digging in the manner intended.

**Indus: From right angled triangles to all triangles (AD1).**

*MK: This proof, let's assume it works. But, then it only works for triangles that can be made from breaking up a rectangle.. into two equal part, right? So, what I was saying to you. So, the question is, is there a way to do it for any triangle? And, she's saying use parallelograms rather than rectangles, right?*

*(A student says something inaudible)*

*MK: Ah..*

*(Other students raise their hands)*

*Arnav: Sir, you can do it for any triangle*

*MK: What is.. your proof? Yeah..*

*Arnav: Consider any triangle*

*MK: Consider any triangle*

*(Student draws a triangle and then drops a perpendicular from one of the vertices, making two*

*right angled triangles – see Figure 16)*

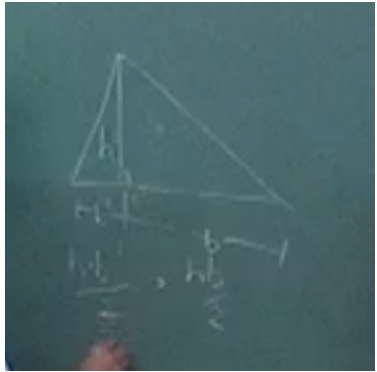

*Figure 16.* Screenshot of green board from conversation with Arnav.

*Arnav: We can drop its altitude. So, this is b.. h. This is b and this is b dash.*

*The base of this triangle is b, and of the other is b'*

*So, the area of this triangle is b dash h by 2 and of this triangle is h b by 2. And for the triangle*

*you get h by two, b dash plus b, which is the base.*

     In this exchange, Arnav gives a proof for the area of a triangle using the area of a right

triangle, an example of AD1. He says that you can drop a perpendicular from the top of a triangle to get two right angled triangles. The calculation of the sum of their area gives us what we want.

This interaction, and what followed, shows the role of the instructor better than almost any other in the two sessions. This was something the instructors had not thought of before the session commenced. They expected the students to demonstrate the area of a triangle via the area of a parallelogram. However, rather than taking the session in the direction they had intended, they let the session go in the direction students took it. This highlights the role of the instructor and requires the instructor to be open to various possibilities and be able to evaluate an argument in real time.

However, the proof doesn't work for all triangles, as I discuss now:

**Indus: Flaw in the proof and a 'fix' (AD1 & AD3).**

*There was a flaw spotted in the proof by Rohit. Arnav's argument works for triangles where the two angles on the base of the triangle are less than 90 degrees. However, if one of them is more than 90 degrees, the perpendicular we drop is outside of the triangle. A student suggested a fix – even though that altitude is outside of the triangle, there always exists at least one altitude which lies within the triangle. This was accepted by the instructors.*

Figure 17 is a representation of Arnav's argument.

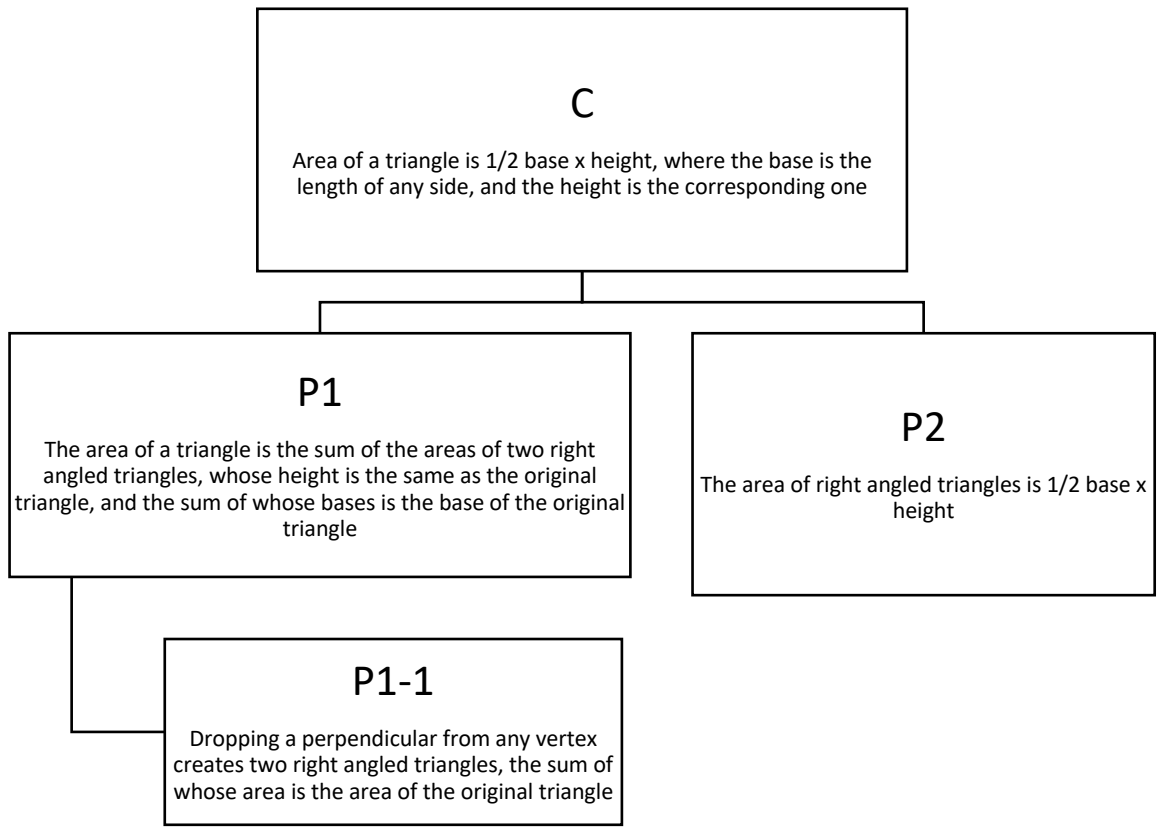

*Figure 17.* Tree representation of Arnav's argument.

For triangles such as that shown in Figure 18, the proof given by Arnav, as represented in Figure 17, doesn't actually work.

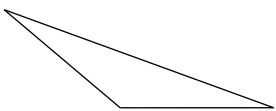

*Figure 18.* Image of counter-example to Arnav's proof.

P1-1 is what the student is questioning, an example of engaging in AD3. The flaw is in the existence of the vertex P1-1, even though P1 does actually follow from P1-1, and C follows from P1 and P2.

However, even the suggested fix, that there will exist at least one vertex whose altitude is inside of the triangle, doesn't solve the problem. The reason it doesn't fix the problem is that in this case we are only concerned with one altitude, the one between the two parallel lines. Hence, we need to show that the area of a triangle is ½ b x h for all three bases and heights, or at least that particular base and height. This flaw was only noticed after the sessions were over when reviewing the footage. So, the instructors accepted an extremely flawed argument and the students were left assuming that the argument works. In this case, there was a conflict between the conceptual and figural aspects, which even the instructors did not notice.

# DISCUSSION

## Engagement with Assumption Digging

The students, at least those of whom spoke during the sessions, clearly valued the Deductive Proof Scheme. Almost never during the sessions did they give a proof in the form of examples or via appeals to authority. This might be a function of the previous workshop sessions or of their experiences before the workshop. Conducting an Assumption Digging session without this would make things a lot more complicated.

The students in the sessions also had already had experience with Euclidean Geometry in school. Nowhere during this session did they come to any new conclusion which they had not seen before. All they were doing was constructing the foundations upon which those conclusions stood.

In both the sessions, students justified the area of a triangle in terms of areas of other shapes. Unlike in the Zaskis and Leikin (2008) study, it did not seem like this was thought of as 'cheating.' Rather, they were giving a proof for a statement starting from something they believed to be true.

The main point of difference between the two sessions, in terms of content, was the proof for the area of a triangle. In Ganga, the session took the expected route – justifying the area of a triangle in terms of the area of a parallelogram, which in turn was justified by the area of a rectangle. In Indus, the route was completely different; parallelograms were not used in the final justification. The proof was in terms of rectangles, and right angled triangles. It could be argued that in Ganga the students were guided to their conclusion by the instructors since their proof was the same as what the instructors had expected. However, at least in Indus, we can conclude that some of the students were thinking for themselves – the instructors had not thought of this proof before the session.

While most of the engagement with Assumption Digging was at the level of extracting premises and definitions on the behest of the instructor, there were times when some students did themselves see the need to ask questions. However, this was restricted to when they found a flaw in the given reasoning rather than when they found a conclusion to be 'obviously' true. A prime example of this was Manu finding a flaw in Bharat's procedure for creating parallelograms out of triangles. Manu was checking to see whether Bharat's procedure actually worked, and found a counter-example (in fact a whole class of counter examples). He was doing what Marrioti and Fiscbein (1997) would call 'returning to observations to check' a definition. Given that this was their first session of Assumption Digging, Manu's skeptical attitude is quite impressive and exactly what this type of session ought to be aiming at.

The other notable event was Tushar's argument for the area of a triangle from the area of a square. What he is doing there is actually laying down a theory, starting by stating his basic assumptions, and then deducing its consequences. While that argument potentially has holes, the remarkable part about it is that it shows an understanding of the structure of a mathematical theory. While this was not made explicit, he was laying out what Marrioti and Fischbein (1997) called the 'basic objects' of the theory which has certain properties. For instance the 'unit' square has the property that its area is 1.

In terms of the graphical representation from the Introduction, most of the AD1 digging was related to vertices with no existing justification. The notable exceptions are where some students found flaws in others' proofs. In those cases, the digging was either on the edges, the steps of reasoning, or on the existence of vertices. Those instances are pointed out in the results section, and the graph diagrams are provided.

## Role of the Instructor

There are various times during the sessions where the instructors played a role in steering the discussions. Many of these have been pointed out in the results section. There were definitely times when the instructors could have asked less leading questions. That could have led to the students doing more of the work. Since the goal of the session was not necessarily to come to a conclusion, even if this would have taken a large chunk of time, it might have been worth it.

However, the instructors also did follow the students' direction at times. An example is the proof for the area of the triangle in Indus. It was clearly important that the instructors were not just following a scripted lesson plan but were willing to go along with what the students suggested.

There were also times when the instructors were wrong. In the first of these instances, with the flipping procedure, this had a valuable result. Manu was able to find the problem and point it out to the entire class. In the second instance, which is at the end of the results section, students were left with a wrong proof. I would contend that this doesn't matter. The goal of the session is not that students have a body of knowledge which is deemed to be correct. Rather, the valuable part of the session is the mindsets and abilities associated with Assumption Digging. What is more important is for instructors to be willing to admit they are wrong and, more importantly, get joy out of being shown to be wrong by their students.

## Classroom Norms

In order to engage in assumption digging, students need to have authority and ought to be held accountable in the sense of Engle & Conant (2002). What is pretty clear from the data is that in this instance, there was a sharing of authority between the instructors and at least some of the students. Students were directly engaging with the mathematics, trying to make sense of it on their own. However, once they put ideas across, they were critiqued, either by their fellow students (as in the case where Tarini rejects the notion of using a square, or in the case where Manu shows that Bharat's flip strategy doesn't work) or by the instructors. These critiques were based on disciplinary norms of mathematical practice as modeled by the instructors in this and previous sessions.

## Instructor Tips

Considering both the data presented in this thesis and my previous experience of teaching similar sessions, the following are some tips which may be useful to mathematics instructors who wish to engage in this work:

### Being Wrong.

As can be seen in the Assumption Digging session described in the thesis, instructors being wrong can have a significant positive impact on a session, especially in getting students to question the instructor. It shifts authority towards students, which is one of the guiding principles laid out by Engle & Conant (2002) for productive disciplinary engagement. While in this case, being wrong was unintentional, this can also be done intentionally. The crucial aspect of this is the instructor's response to being wrong. The way I tend to deal with it is to admit I'm wrong and celebrate that students were able to demonstrate that. However, this needs to be genuine. In my experience, people can see through faked appreciation and humbleness.

### Exploring New Math.

When walking into the sessions, the instructors did not have a well worked out plan or even a clearly articulated proof for the claim. We were coming up with things along with the students. Not having a proof made it easier for us to not push students down a particular path. It placed the instructors on a much more equal footing to that of the students than it would have if we had entered the classroom with a clearly worked out proof. It also resulted in us being wrong naturally. This does not imply no preparation is required. It is just that the preparation required is not in terms of the mathematical content of the session, but rather in terms of getting an understanding of the structure of mathematics and of mathematical theories.

**Attempt to ask questions which are less leading.**

When students appeared to be at an impasse, one commonly adopted strategy by the instructors was to ask them a question. At times during these sessions, as noted above, the instructors asked questions which directly led students to the answers. In some cases, this could result in a missed learning opportunity. Rather, instructors should attempt to ask questions which push students in a fruitful direction without constraining their options more than necessary.

**Euclidean Geometry**

These Assumption Digging sessions highlight some of the failures of the regular Euclidean Geometry course, as documented by Weiss and Herbst (2015). The students had clearly accepted many claims, such as conclusions about area, without proof. They did not have a clear definition of equidistant before the workshop even though they would have used that notion in their Euclidean Geometry course. This session at least attempts to achieve the stated goals of Euclidean Geometry courses. Assumption digging sessions, using what students have learnt in regular Euclidean Geometry courses as a starting point, have the potential to give students an understanding and experience of doing real mathematics.

There were some indications of where some students were in relation to the Van Hiele levels. The example of Tushar, expanded upon above, demonstrates some understanding of axioms, definitions, theorems and proofs, which corresponds to level 3. It also shows some understanding of the relationship between axioms and theorems, which is a part of level 4.

**Limitations of the Study and Future Work**

The students who participated in these sessions already seemed to see the value in deductive proofs. This may not be the norm amongst students more generally (Harel & Sowder, 2007; Weber, Inglis, & Mejia-Ramos, 2014). Hence, the sessions may look very different if the lesson plan is used elsewhere without giving students other related experiences.

While the mindsets and abilities associated with Assumption Digging are mentioned in this thesis, these definitely need unpacking and elaboration in order to better evaluate an Assumption Digging session.

The instructors in these sessions are a math graduate student and a professional theoretical linguist. Hence, both have a lot of experience evaluating mathematical arguments. Many teachers may not have had such experiences, and hence it may be hard for them to engage skeptically with students' arguments. So, it is not clear as to whether such a session could be a part of the regular school curriculum without some effort being made to work with teachers.

Finally, this was only one session in each school. It would be interesting to see how thing develop as students have more and more experience with Assumption Digging.

It is important that such sessions be tried out in other contexts and with other instructors in order to say whether Assumption Digging ought to be a part of the curriculum. Of course, one session will not be sufficient, and it will probably have to form just a part of a larger course on Theory Building in Mathematics.

**Relationship with Areas outside of Math**

Unpacking arguments outside of mathematics, while not identical, involves similar strategies. Suppose we are evaluating the following claim made by a student: 'Athena is more intelligent than me because she did better on an exam.' A premise here is that 'those who do better in this exam are more intelligent.' In order to evaluate whether we should accept that claim, we need clarity on what we mean by intelligence. Then, we need to show that this exam tests for intelligence. Finally, we need to evaluate the accuracy of the exam, asking questions like how much scores vary over successive takes, etc. The forms of reasoning used here are not the same as in mathematics since we allow, for instance, sample to population reasoning. However, the broader strategy is similar. It is unlikely to be the case that just doing Assumption Digging in math will have a huge impact outside of math. However, it is plausible that doing Assumption Digging in math and outside, and talking about the connection, could have an even greater impact on those areas. Of course, this is just plausible speculation.

**Other Forms of Theory Building**

Some of the other forms of Theory Building mentioned in the introduction were amongst the other sessions of the same workshops. While the main difference in the types is in their motivation, there are significant similarities in their practice. You will notice that in the lesson plans in Appendix 2, and in Fawcett (1938), there are instances of Assumption Digging occurring during sessions where other aspects of Theory Building are being highlighted. The next step in exploring theory building in education would be to construct and try out an entire course on the topic, which would include all the different types of theory building. A course over a longer period would also give us a better idea as to whether there are any changes in students' mindsets which are visible to those who interact with them outside the class, such as their other teachers and parents.

## APPENDIX 1 – PARALLEL LINES

The following is a lesson plan for a session on defining Parallel Lines. While this is indicative of what happened in one of the pre-Assumption Digging sessions in the workshop, the actual session diverted significantly from it. The lesson plan is written as a dialogue between a teacher and a student.

Teacher: What are parallel lines?

Student: They are lines which do not meet.

Teacher: Are these lines parallel?

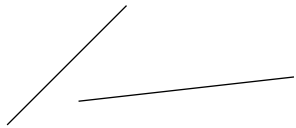

Student: No. They are segments, not lines. Lines go on forever.

Teacher: So are you saying segments cannot be parallel. Are the two segments below parallel?

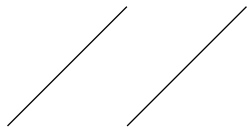

Student: Yes they are. How about if we say that 2 lines or segments are parallel if they do not meet even when extended.

Teacher: Okay. How about the following. Are they parallel lines?

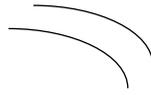

Student: No, they are curves, not lines.

Teacher: I'm guessing that when you say line, you mean straight path and when you say curve, you mean a non-straight path. Is that right?

Student: Yes.

Teacher: The words are a little confusing since we seem to have two commonly used terms for the same thing: line and straight line. Just for the purposes of this session, let us use the following classification:

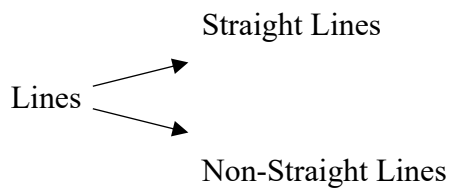

A line is something you can trace with your finger without lifting it (of course in the geometry we are working in, you cannot actually trace a line since it has no width). It could be finite or not. We can call segments finite lines.

Student: Does that mean that when you ask us for what parallel lines are, you want us to include non-straight lines?

Teacher: That is your decision. We would like the definition to be as general as possible. So, if you can come up with aa definition which covers both straight and non-straight lines and leads to interesting theorems, then you should do that. If you think concentrating on straight lines is better, you can do that. In order to make that decision, let us put down some common alternative definitions of parallel lines along with some scenarios:

Definition 1: Two lines are parallel if they do not meet.

Definition 1': Two lines are parallel if they do not meet when extended.

Definition 2: Two lines are parallel if they are equidistant.

Definition 3: Two lines are parallel if there exists a third line which is perpendicular to both.

Student: How do we evaluate these.

Teacher: Notice that Definition 1 is useless given the scenario we already put forward above. By Definition 1, any two paths would be parallel as long as they don't share a point in common.

Student: What about Definition 1'. What does it mean to extend a non-straight line? See the following. How would we extend it?

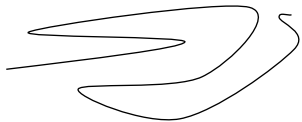

Teacher: Great observation. Unless we can come up with some way of extending any line, Definition 1' only makes sense for straight lines since it is obvious how we can extend straight lines.

Student: What about definition 2?

Teacher: What exactly do we mean by equidistant? This is quite a confusing notion for non-straight line, so let us put it aside for now. Let's try Definition 3. To illustrate definition 3, let us see what it means for straight lines:

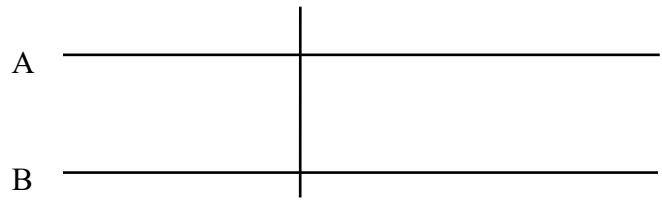

Student: This seems to work for straight lines. What does it mean for a line to be perpendicular to a curved line.

Teacher: Does it work for straight lines? Let me give you a scenario:

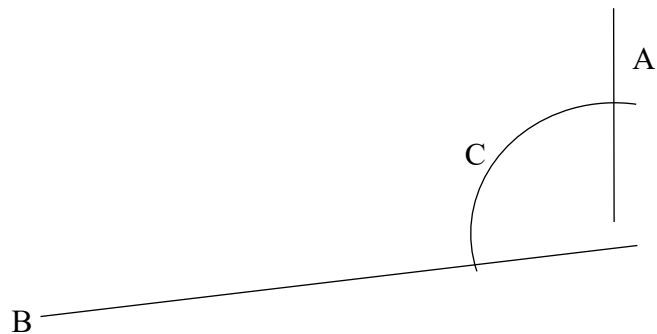

Both A and B are perpendicular to C. However, would you judge A and B to be parallel.

Student: No. But, we can fix the definition to say:

Definition 3': Two lines are parallel if there exists a **straight** line which is perpendicular to both.

Teacher: Let me illustrate a consequence of that:

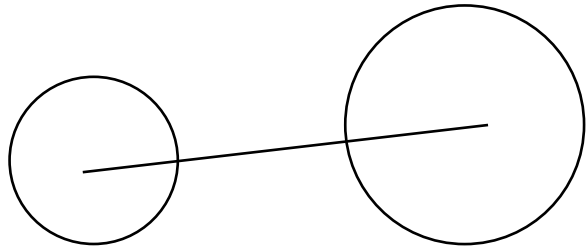

If we draw a straight line connecting the center of any two circles, it is perpendicular to both. Hence, any two circles are perpendicular.

Student: How is the line connecting the centers perpendicular to the circle? What does it mean for a straight line to be perpendicular to a circle?

Teacher: Good point. In order to understand this definition, we need to understand what it means for two lines to be perpendicular. To start off with, what does it mean for two straight lines to be perpendicular?

Student: They are perpendicular if there is a 90 degree angle between them.

Teacher: 90 degrees is just a matter of a measurement tool. You can measure lines in the sort of geometry we are doing, right? You probably heard that lines have no width, so how can you measure them? Also, Napoleon used a 100 degree measure rather than a 360 degree measure. If he had won, we would probably be using that now. In that case, what we now call 90 degrees would be 25 degrees. Can you come up with a definition which doesn't require us to measure?

Student: Can't think of one.

Teacher: Let me give one and you can see if it works. We start with saying that two lines are perpendicular if they are at right angles. Now, we have to define right angle. We say a right angle is a quarter of a full rotation.

Student: That makes sense. In the measurement as well, 90 degrees is a quarter of 360 degrees. However, it still doesn't make sense for the circle stuff.

Teacher: How about the following: Line A is perpendicular to line B if the portion of line B on one side of line A is a reflection of the portion on the other side. By that, what I mean is put a mirror on line A. If what you get when you look through the mirror is exactly line B, then line A is perpendicular to line B.

Student: What happens if A is curved?

Teacher: That's interesting. It doesn't seem to make sense if A is curved. So, for now let us say that A has to be straight. Later, you might want to explore a scenario where A is curved.

Student: So, a consequence of that is the following are perpendicular:

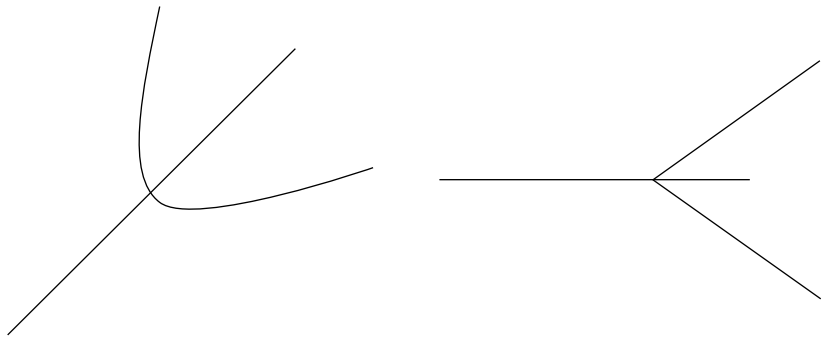

Teacher: You're right. Let us accept this for now. If we get back to the circles, are you okay with accepting that they are parallel?

Student: It seems quite weird.

Teacher: It is worth exploring parallelness of curves more carefully. However, for now, let us stick to straight lines. Before we move on, as an interesting aside, notice that our definition of perpendicular allows for line A to be perpendicular to line B while line B isn't perpendicular to line A. For straight lines, it seems obvious that A perpendicular to B implies B perpendicular to A. Can you come up with a proof for that?

Moving on, so far we have been looking at parallel lines on flat surfaces. Let us move from flat surfaces to curved surfaces, specifically spheres. What do you think parallel lines on spheres look like? Think about globes of the Earth you may have seen.

Student: Well, latitudes look like they are parallel. They are equidistant from each other. They also do not meet and I think any longitude would be parallel to them by our definition.

Teacher: How about if we restrict to straight lines on a sphere? Are latitudes straight lines?

Student: No, they are circles.

Teacher: They might be circles, but are they straight lines?

Student: How can a circle be a straight line?

Teacher: What is a straight line?

Student: What we learnt in class was that a straight line is the shortest path.

Teacher: Take a model of the globe and take a string. Pick any two points on the same latitude (pick them in the same hemisphere but not too close and not the equator). Make sure the string is going through the two points and tighten the string so that it gets as short as it can be. You will notice that the string does not go through the latitude. However, if you try longitudes, you will see that they are the shortest paths. In fact, not just longitudes but any great circle is a straight line on a sphere. A great circle is a line which connects any two antipodal (opposite) points on a sphere. So, what is the consequence if we restrict our definitions of parallel lines to these 'straight lines'. Are there parallel lines?

Student: Well, all great circles meet so by that definition they are not parallel. They are not equidistant either so they fail by that definition as well.

Teacher: However, if you take the third definition, then any pair of great circles are parallel. To see that look at them drawn on a globe. You will notice that they intersect at two antipodal points. Think of these as the 'North' and 'South' Poles. Then there is a great circle which is the equivalent of the equator. We can draw that on the globe. You can see that it is perpendicular to both the great circles we were interested in.

Student: So, the definitions resulted in the same consequences for straight lines on a flat surface but not on a sphere.

Teacher: Exactly! That is why you need to be careful with checking definitions. That is also why you should always try to use a definition in areas where it originally wasn't intended – it can have interesting consequences. I will leave you with a few questions:

1. On a flat surface, the sum of angles of a triangle is two right angles (180 degrees). What is the sum of the angles of a triangle on a sphere? Remember a triangle has to be surrounded by straight lines. Think about other shapes as well.

2. You had remarked earlier that 'circles can't be straight lines.' Great circles are straight lines by your own definition. The interesting question is 'are they circles?' If they are, then circles can be straight lines. What about latitudes? Are they circles on a sphere? In order to answer that you need to define what a circle is. Remember that when we are working on a sphere, we want to stay within the sphere. So, the center of the circle cannot be inside the sphere. It must be on the surface.

## APPENDIX 2 – LESSON PLANS FROM OTHER WORKSHOP SESSIONS

Appendix 1 contains the lesson plan for the session on Parallel Lines. In this appendix, there are lesson plans for a few other sessions from the same workshop. The parallel lines session has been separated from the rest since there is direct reference to it in the results section. The reason for including these other lesson plans is to give the reader a sense of what the students were exposed to before the Assumption Digging session.

**Straight Lines and Intersections**

### Part 1

Teacher: We will be working on a flat surface. What you have in this world are points and straight lines. We have two possibilities of worlds based on their size:

1.  a finite world with definite boundaries, like a sheet of paper.

2.  an infinite flat world which goes on forever in all directions

We also have three possibilities, based on the length of the straight lines

A.  only finite straight lines allowed (what are also called line segments)

B.  only infinite straight lines allowed, which go on forever in both directions

C.  both finite and infinite lines allowed

If not specified later, by default, assume we are working in a world which satisfies condition 2 and condition C. These are some of the axioms of this world. An axiom is a statement we assume to be true and work from there.

I am setting one other condition in this world: not more than two lines can intersect at a given point. This gives us another axiom. Assuming the conditions set so far, let me ask you a questions: Given 4575 straight lines and exactly 25 points of intersection per line, what are the total number of points of intersection?

Student: This seems like a very hard problem. It will take too long to figure it out.

Teacher: What a mathematician would do when confronted with such a problem, would be to first generalize it. What that means is that rather than solving a problem about 4575 lines and 25 points of intersection per line, we try to solve the problem about $n$ lines and $i$ points of intersection per line.

Student: That looks nicer, but seems even more difficult. At least with the large number of lines, we could spend a long time and eventually come up with  solution.

Teacher: Well, what mathematicians are looking for are general patterns, such that our specific problem becomes much easier to solve. At this stage, a mathematician would try out simple examples of $n$ lines and $i$ points of intersection per line. All we can do is to hope that there is some general pattern which will help us with the specific example at hand.

Student: So, you mean trying out values of $n$ and $i$ like 2 and 1 and so on?

Teacher: Yes. Lets put down some examples to explore:

| S. No. | $n$ | $i$ |
|--------|-----|-----|
| 1      | 1   | 0   |
| 2      | 1   | 1   |

What are the total number of points of intersection here?

Student: The first one has zero points of intersection.

Teacher: How do you know that?

Student: Well there is only one line, and that line has zero points of intersection, so the total number of points of intersection must be zero, right?

Teacher: Sure, so lets move on to the second example.

Student: Hmm.. I don't think the second example is possible to make.

Teacher: Why is that?

Student: Since there is only one straight line, and straight lines cannot intersect with themselves, one straight line can have at most zero points of intersection.

Teacher: How do you know that straight lines cannot intersect with themselves?

Student: Isn't that obvious?

Teacher: In mathematics, the word obvious is not allowed. You have to give reasons for everything.

Student: I can't think of any reason, so can't we just let this be one of those axiom things?

Teacher: Yes we can, for now. However, spend some time thinking about how you can define straight lines. Lets leave that for a different lesson.

Student: So, it seems like there might be some examples where it is not even possible to create the configuration - forget about counting the number of points of intersection. Maybe the example you gave initially is of the same sort?

Teacher: Maybe. So, it might be time to change the question we are asking. As we have seen, we should first be asking a related but similar question: For what values of *n* and *i* is it possible to create a configuration of *n* lines with exactly *i* points of intersections per line? How about *n=1* and *i=2*?

Student: That is not possible. We already showed it earlier when we said that one line can have most zero points of intersection per line.

Teacher: Great! What you are quoting from earlier is a theorem. Let's state it clearly:

*Theorem 1: Given any configuration with exactly one straight line, there are zero points of intersection in this configuration.*

Student: So, can we state this new result also as a theorem?

Teacher: Yes. A theorem is any statement you have proved from the definitions and axioms.

*Theorem 2: If n=1 and i is greater than or equal to 1, then a configuration with n lines and i points of intersection per line is impossible to create.*

Student: So, lets move on to more examples. How about *n=2* and *i =1*. That one clearly is possible. Just take two straight lines and make them intersect.

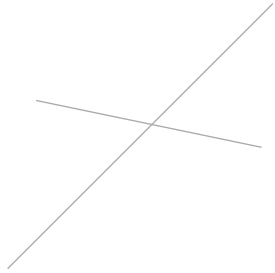

Teacher: Great!

*Theorem 3: A configuration with 2 lines and 1 point of intersection per line is possible to create.*

Student: However, 2 lines and 2 points of intersection per line seems to be impossible.

Teacher: For now, that counts as a conjecture. A conjecture is a statement without a proof. When you prove a conjecture, you get a theorem.

Student: Well, start with a line. That line cannot intersect with itself from, which we get from Theorem 1. So, it has to intersect with two other lines. There is only one more line allowed, so this configuration is impossible.

Teacher: Why does it have to intersect with two other lines? Why can't it intersect with one other line twice?

Student: Thats impossible. But, I can't prove it. So, lets make it another axiom?

Teacher: Sure!

*Theorem 4: If n=2 and i is equal to 2, then a configuration with n lines and i points of intersection per line is impossible to create.*

Student: I'm beginning to have a suspicion here that *n* has to be greater than *i.*

Teacher: Another conjecture! Can you prove it?

Student: I think so. From our axioms, we already know that one line cannot intersect with itself. Also, two lines can have at most one point of intersection. Given *n* lines, pick a line. There are only *n-1* lines left. So, the line you picked can have at most *n-1* points of intersections.

Teacher: Fantastic! Do you realize, you have done what seems to be infinite work in a finite amount of time! You have not only proved something for *n=2* and *i=2*, but also for *n=321* and *i*

$=342423$, and an infinite number of other cases. One of the greatest parts of being a

mathematician is that your goal is to do the most amount of work with the least amount of effort.

*Theorem 5: If i ≥ n, then a configuration with n lines and i points of intersection per line is*

*impossible to create.*

Student: So, lets explore some more examples.

Teacher: Your homework is to try out the following:

| S. No. | n | i | Possible? | S. No. | n | i | Possible? |
|---|---|---|---|---|---|---|---|
| 3 | 3 | 0 | | 13 | 5 | 3 | |
| 4 | 3 | 1 | | 14 | 5 | 4 | |
| 5 | 3 | 2 | | 15 | 6 | 0 | |
| 6 | 4 | 0 | | 16 | 6 | 1 | |
| 7 | 4 | 1 | | 17 | 6 | 2 | |
| 8 | 4 | 2 | | 18 | 6 | 3 | |
| 9 | 4 | 3 | | 19 | 6 | 4 | |
| 10 | 5 | 0 | | 20 | 6 | 5 | |
| 11 | 5 | 1 | | | | | |
| 12 | 5 | 2 | | | | | |

**Part 2**

Teacher: Did you figure the homework out?

Student: A few of them. I think I can state a conjecture: A configuration with *n* lines and zero points of intersection per line is always possible to create.

Teacher: How would you go about constructing such a configuration?

Student: Just take *n* finite straight lines which do not intersect. You can do it with infinite lines as well, but this is easier.

Teacher: Great.. Another theorem

*Theorem 6: A configuration with n lines and zero points of intersection per line is always possible to create.*

Student: I don't think I have any more general ones.

Teacher: Thats fine. Lets go through the list one by one. Maybe, we will see a pattern. Before we do that, lets introduce some notation. Lets say *(n,i)* represents a configuration with *n* lines and *i* points of intersection per line. So, lets start with *(3,1)*.

Student: I haven't been able to create it. I think its impossible, but don't have a proof.

Teacher: Lets come back to that later. How about *(3,2)*.

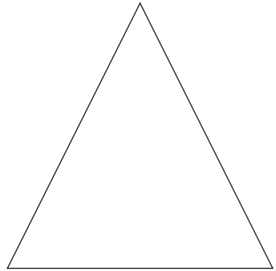

Student: Thats possible. Just draw a triangle.

Teacher: *Theorem 7: (3,2) is possible.*

Student: *(4,1)* is also possible

Teacher: How?

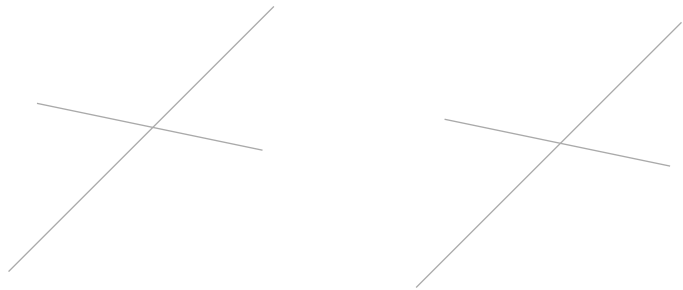

Student: By drawing two crosses:

Teacher: *Theorem 8: (4,1) is possible to create*

If we restricted ourselves to only infinite lines, would it still be possible?

Student: I don't think so

Teacher: I'm going to leave that to you to figure out why

Student: *(4,2)* is also possible by making a square, and I think I've found a pattern. *(n,2)* is always possible for all *n*>2

Teacher: How would you make *(n,2)* for an arbitrary *n*

Student: By making an *n*-sided polygon

Teacher: *Theorem 9: (n,2) is possible to create for all n > 2*

Student: *(4,3)* is also possible:

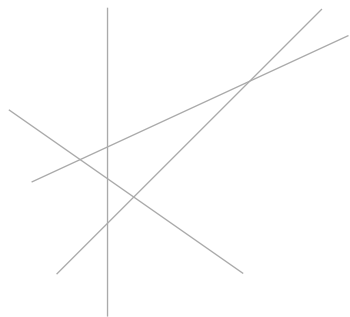

Teacher: *Theorem 10: (4,3) is possible to create*

Student: That gives me an idea - I think *(n,n-1)* is always possible for all *n*

Teacher: Can you prove that?

Student: Can't you just take *n* infinite lines? Each one will intersect with all the others resulting in *n-1* points of intersection.

Teacher: Can you pick any *n* lines? Do two infinite lines always intersect with each other?

Student: Oh! Not if they are parallel. So, take *n* non-parallel lines.

Teacher: For that you need to define parallel lines.

Student: They are lines which do not intersect.

Teacher: Let's see what you are saying:

I asked whether infinite lines always intersect. You replied that they do not do so if they are parallel. Then when I asked you what parallel lines are you replied that they are lines which do not intersect. The reasoning here seems circular, doesn't it?

Student: Then, can't we say parallel lines are those lines which are equidistant.

Teacher: In that case, what is distance? We don't have a way to measure that in this world

Student: Then how do we deal with this?

Teacher: Mathematicians have developed sophisticated ways to introduce the concept of distance to such worlds. However, we can play a trick Euclid, the first real mathematician we know of, often played. We can say parallel lines are undefined entities (what are also called primitives) with the following properties (which are also axioms of the world):

    a)   given any straight line, there are infinite parallel lines to it

    b)   parallel lines do not intersect with each other

    c)   infinite non-parallel lines always intersect with each other

Student: We will also need to add in another axiom: given a straight line, there are infinite non-parallel lines to it. That should give us a proof that *(n,n-1)* is possible for all n.

Teacher: There is still one problem. Remember that our world requires at most two straight lines to intersect at a given point. We need to show that it is possible to avoid this when creating the *(n,n-1)* configuration. Lets leave that aside for now and let *(n,n-1)* be a theorem contingent on us proving this intermediate result:

*Theorem 11: (n,n-1) is possible to create for all n > 0*

Student: So far, apart from when $i \geq n$, we have proved everything else so far as possible. Maybe *(n,i)* is always possible if i < n.

Teacher: Thats another conjecture. However, is it true? Earlier, you said *(3,1)* might not be possible to create. Can you create it?

Student: No. Not so far.. Actually, I don't think it's possible

Teacher: Can you prove that?

Student: Well, take two lines and make them intersect. You have one left over. That line has nothing to intersect with. So, *(3,1)* is not possible.

Teacher: I don't think that reasoning works. What you need to prove is that no matter what you do, *(3,1)* is impossible to create. You have given one method and shown it is impossible using that method. Maybe, there is some other way to accomplish it.

Student: I don't think there is any way. How does this reasoning sound:

Take a line. We need it to have one point of intersection. Since it cannot intersect with itself, it will have to intersect with one of the other two lines. Make that happen. Now, two lines have exactly one point of intersection. We know that since we have an axiom that two lines can intersect at most at one point. We are left with one line. That line cannot intersect with any of the other two lines since if it did, at least on line would have two points of intersection. So, there is no way for that line to have more than zero points of intersection. Hence, *(3,1)* is not possible to create.

Teacher: Great! As you might have noticed, in mathematics, proving something to not be possible is usually much harder than proving that it is possible.

Student: I think I have found another pattern - *(n,1)* is possible if *n* is even and *(n,1)* is impossible if *n* is odd.

Teacher: Let me introduce some terminology. When mathematicians state something like you have, they usually use the phrase 'if and only if.' So, your claim would become:

*(n,1)* is possible to create if and only if *n* is even. What this means is two things:

    *1.* *(n,1)* is possible to create if *n* is even

    *2.* If *(n,1)* is possible to create, then *n* is not even

This is the same as what you said. Can you come up with a proof?

Student: To show that *(n,1)* is possible if *n* is even, you can just create *n/2* crosses with finite lines which do not intersect with each other. To show *(n,1)* is impossible if n is odd, we will have to show that you are forced to create crosses till you have just one line left.

## Discrete Geometry

### Part 1

This session starts with an ill-defined question:

Given a world with exactly 6 points, can every straight line be bisected?

Why is this question ill-defined? It is ill-defined since we don't know what we mean by a world with exactly 6 points, we don't know what a straight line in that world looks like, or if they even exist, and we don't have a definition for bisected which is directly transferrable to this world. Hidden behind both the notions of straight line and bisection is the concept of distance, which is also not defined.

The goal of this portion of the session is for students to see the relationship between definitions, assumptions, and conclusions. Their task is to define these concepts, and follow the consequences of the definitions in order to come to a conclusion about the question.

The following are some potentially helpful suggestions for students if they need them:

### Concept of Neighbor

In Euclidean Geometry, points don't have neighbors. This is because between any two points, there exist more points. However, in this case, given exactly 6 points, this concept might be useful.

### Suggested Starting Configuration

While there are many configurations of neighbor relations that students could choose from given 6 points, a suggestion could be that they start with a simple one – each point except from two of them, who have one neighbor, have exactly two neighbors (points are arranged in a 'line')

**Distance**

It may not be clear from the question that the notion of distance is crucial. A starting point for defining distance could be that the distance between any two neighbors is the same, say 1 (they could also think about these distances being different, but this should probably be left for later in the session). Then, the distance between two points will be the number of neighbors in the shortest (in terms of count) sequence of neighbors connecting those two points.

**Straight Line**

Straight line between two points could be defined as the shortest path of points connecting those two points.

**Representations**

If we are dealing with the simple situation of a distance of 1 between two neighbors, the representation of a graph might be a useful thing to introduce. Here, vertices would describe the points, and edges would signify the relationship of neighbor. The 'length' of the edges in the drawing of the graph bears no relationship to the notion of distance.

**Bisection**

There are two ways to define bisection in a way that it bears resemblance to the notion of bisection from Euclidean Geometry. The first is breaking a straight line into two straight lines of equal length. The other requires a straight line which is the bisector. These two definitions will have distinct consequences.

**Part 2**

In this section, we explore various possible worlds we can create using some of the notions we defined above in Part 1. We can generalize to worlds with n points rather than just 6, look at various configurations of these n points, and potentially allow for neighbors to have different distances between them. While students should be free to explore these consequences, it ought to be suggested that they don't make the worlds overly complicated too soon. The following are questions we can inquire into:

**Circles**

Given the definition of circles from Euclidean Geometry, what do circles in these worlds look like? Are there worlds of this types you can create where circles do not exist?

**Polygons**

What sorts of polygons exist in these worlds (given the notion of polygon transferred from Euclidean Geometry)? Try to explore what regular polygons look like. Can you create worlds where certain types of polygons do not exist?

**REFERENCES**

Bass, H. (2005). Mathematics, mathematicians, and mathematics education. *Bulletin of the*

      *American Mathematical Society, 42(04),* 417–431. https://doi.org/10.1090/S0273-0979-

      05-01072-4

Bass, H. (2017). Designing opportunities to learn mathematics theory-building practices.

      *Educational Studies in Mathematics, 95(3),* 229–244. https://doi.org/10.1007/s10649-

      016-9747-y

Bolyai and Lobachevsky. (n.d.). Retrieved from

      https://www.storyofmathematics.com/19th_bolyai.html.

Euler. (n.d.). Retrieved from https://www.storyofmathematics.com/18th_euler.html.

Christofferson, H. C. (1930). A fallacy in geometry reasoning. *The Mathematics Teacher*, *23*(1),

      19-22.

Coppin, C. A., & Mathematical Association of America (Eds.). (2009). *The Moore method: A*

      *pathway to learner-centered instruction*. Washington, DC: Mathematical Association of

      America.

Dawkins, P. C. (2015). Explication as a lens for the formalization of mathematical theory

      through guided reinvention. *The Journal of Mathematical Behavior*, *37*, 63–82.

      https://doi.org/10.1016/j.jmathb.2014.11.002

Fawcett, H. P. (1938). The nature of proof: A description and evaluation of certain procedures

      used in senior high school to develop an understanding of the nature of proof. *National*

      *Council of Teachers of Mathematics Yearbook*, *13*, 1-153.

Flener, F. O. (2009). The "guinea pigs" after 60 years. *Phi Delta Kappan*, *91*(1), 84–87.

      https://doi.org/10.1177/003172170909100120

González, G., & Herbst, P. G. (2006). Competing arguments for the geometry course.

    *International Journal for the History of Mathematics Education, 1(1),* 7-33.

Gowers, W. T. (2000). The two cultures of mathematics. In V. Arnold et al. (Eds.), *Mathematics:*

    *Frontiers and Perspectives* (pp. 65–78).

Harel, G., & Sowder, L. (2007). Toward comprehensive perspectives on the learning and

    teaching of proof. In F. Lester (Ed.), *Second Handbook of Research on Mathematics*

    *Teaching and Learning*, 805–842.

Hilbert, D. (2014). Foundations of geometry. Retrieved from

    http://www.myilibrary.com?id=802854

Larsen, S. P. (2013). A local instructional theory for the guided reinvention of the group and

    isomorphism concepts. *The Journal of Mathematical Behavior*, *32*(4), 712–725.

    https://doi.org/10.1016/j.jmathb.2013.04.006

Lockhart, P. (2009). A mathematician's lament: How school cheats us out of our most

    fascinating and imaginative art form. *New York, NY: Bellevue Literary Review.*

Mariotti, M., & Fischbein, E. (1997). Defining in classroom activities. *Educational Studies in*

    *Mathematics*, *34*(3), 219–248.

National Council of Educational Research and Training (India). (2006). *National curriculum*

    *framework 2005: position paper.* New Delhi: National Council of Educational Research

    and Training.

Polya, G. (1963). On learning, teaching, and learning teaching. *The American Mathematical*

    *Monthly*, *70*(6), 605. https://doi.org/10.2307/2311629

Polya, G. (1978). Guessing and proving. *The Two-Year College Mathematics Journal*, *9*(1), 21.

    https://doi.org/10.2307/3026553

Polya, G. (1979). More on guessing and proving. *The Two-Year College Mathematics Journal*, *10*(4), 255. https://doi.org/10.2307/3026620

Rao, S. (1917). *T. Sundara Row's geometric exercises in paper folding.* The Open Court Publishing Co, 148.

Riemann. (n.d.). Retrieved from https://www.storyofmathematics.com/19th_riemann.html.

Russell, B. (1919). The study of mathematics. *Mysticism and Logic and Other Essays*. London: George Allen & Unwin Ltd.

Scheiner, T. (2016). New light on old horizon: Constructing mathematical concepts, underlying abstraction processes, and sense making strategies. *Educational Studies in Mathematics*, *91*(2), 165–183. https://doi.org/10.1007/s10649-015-9665-4

Standards for Mathematical Practice. (n.d.). Retrieved from http://www.corestandards.org/Math/Practice/

Tall, D. (2008). The transition to formal thinking in mathematics. *Mathematics Education Research Journal*, *20*(2), 5–24. https://doi.org/10.1007/BF03217474

Tall, D., & Vinner, S. (1981). Concept image and concept definition in mathematics with particular reference to limits and continuity. *Educational Studies in Mathematics*, *12*(2), 151–169. https://doi.org/10.1007/BF00305619

Tall, D., Yevdokimov, O., Koichu, B., Whiteley, W., Kondratieva, M., & Cheng, Y.-H. (2011). Cognitive development of proof. In G. Hanna & M. de Villiers (Eds.), *Proof and Proving in Mathematics Education* (Vol. 15, pp. 13–49). https://doi.org/10.1007/978-94-007-2129-6_2

Usiskin, Z. (1980). What should not be in the algebra and geometry curricula of average college-bound students? *The Mathematics Teacher*, *73*(6), 13.

Van den Heuvel-Panhuizen, M., & Drijvers, P. (2014). Realistic mathematics education. In S.

Lerman (Ed.), *Encyclopedia of Mathematics Education* (pp. 521–525).

https://doi.org/10.1007/978-94-007-4978-8_170

Vojkuvkova, I. (2012). The van Hiele model of geometric thinking. WDS'12 Proceedings of

Contributed Papers, Part I, 72–75

Weber, K., & Alcock, L. (2004). Semantic and syntactic proof productions. *Educational Studies*

*in Mathematics*, *56*(2/3), 209–234.

Weber, K., Inglis, M., & Mejia-Ramos, J. P. (2014). How mathematicians obtain conviction:

Implications for mathematics instruction and research on epistemic cognition.

*Educational Psychologist*, *49*(1), 36–58. https://doi.org/10.1080/00461520.2013.865527

Weiss, M., & Herbst, P. (2015). The role of theory building in the teaching of secondary geometry. *Educational Studies in Mathematics*, *89*(2), 205–229.

https://doi.org/10.1007/s10649-015-9599-x

Weiss, M., Herbst, P., & Chen, C. (2009). Teachers' perspectives on "authentic mathematics"

and the two-column proof form. *Educational Studies in Mathematics*, *70*(3), 275–293.

https://doi.org/10.1007/s10649-008-9144-2

Whitehead, A. N., & Russell, B. (1912). *Principia mathematica (Vol. 2).* University Press.

Wubbels, T., Korthagen, F., & Broekman, H. (1997). Preparing teachers for realistic

mathematics education. *Educational Studies in Mathematics*, *32(1),* 1-28.

Zandieh, M., & Rasmussen, C. (2010). Defining as a mathematical activity: A framework for

characterizing progress from informal to more formal ways of reasoning. *The Journal of*

*Mathematical Behavior*, *29*(2), 57–75. https://doi.org/10.1016/j.jmathb.2010.01.001

Zazkis, R., & Leikin, R. (2008). Exemplifying definitions: a case of a square. *Educational Studies in Mathematics*, *69*(2), 131–148. https://doi.org/10.1007/s10649-008-9131-7